\documentclass[12pt]{article}
\usepackage[utf8]{inputenc}
\usepackage[T2A]{fontenc}
\usepackage[russian,english]{babel}
\usepackage{amsmath,amsfonts,amssymb,euscript,graphicx,wrapfig,multirow}
\usepackage{dsfont}
\usepackage{amsthm}
\inputencoding{utf8}
\usepackage{cite}

\usepackage{hyperref}

\makeatletter
\renewcommand{\@biblabel}[1]{#1.}
\makeatother

\theoremstyle{theorem}
\newtheorem{theorem}{Theorem}
\theoremstyle{dfn}
\newtheorem{dfn}{Definition}
\theoremstyle{remark}
\newtheorem{remark}{Remark}

\begin{document}

\selectlanguage{english}

\title{
	On particular diameter bounds for integral point sets in higher dimensions
	\footnote{
		This work was carried out at Voronezh State University and supported by Russian Science
		Foundation grant 19-11-00197.
	}
}

\author{
	Avdeev N.N.
	\footnote{nickkolok@mail.ru, avdeev@math.vsu.ru}
	, Zvolinsky R.E., Momot E.A.
	\\
	\\
	\emph{Voronezh State University}
}

\maketitle

\paragraph{Absctract.}
	An integral point set $\mathcal{P}$ is a set of $n$ points in
	the $m$-dimensional Euclidean space $\mathbb{R}^{m}$ with pairwise
	integral distances,
	such that $\mathcal{P}$ is not contained in a $(m-1)$-dimensional hyperplane.
	In the present paper we discuss some classes of planar integral point sets ($m=2$);
	we mostly focus on the subsets of a union of two straight lines:
	facher, rails, scissors and pyramid sets.
	We use rails and pyramid sets to construct lower bounds for a diameter
	of integral point sets in higher dimensions.

\section{Introduction}
\label{sec:intro}

The study of integral point sets belongs to the classical analysis.
It is motivated by both pure theoretical aspects~\cite{ascher2020erdos,solymosi2010question}
and applications~\cite{kurz2007enumeration}.

\begin{dfn}\label{dfn1}
	A planar integral point set (PIPS) is a set $\mathcal{P}$
	of non-collinear points in the plane $\mathbb{R}^{2}$ such that
	for any pair of points $P_{1}, P_{2} \in \mathcal{P}$
	the Euclidean distance $|P_{1}P_{2}|$
	between points $P_{1}$ and $P_{2}$ is integral.
\end{dfn}

Definition \ref{dfn1} can be generalized as follows:

\begin{dfn}
	An integral point set $\mathcal{P}$ is a set of $n$ points in
	the $m$-dimensional Euclidean space $\mathbb{R}^{m}$ with pairwise
	integral distances,
	such that $\mathcal{P}$ is not contained in a $(m-1)$-dimensional hyperplane.
\end{dfn}

How can we characterize an integral point set?
First of all, we can look at its dimension;
then, we can find its cardinality, which is always finite~\cite{anning1945integral,erdos1945integral};
finally, we can compute the diameter of a finite point set,
which is naturally defined as follows:

\begin{dfn}
	The diameter of an integral point set $\mathcal{P}$ is defined as
	\begin{equation}
		\operatorname{diam(\mathcal{P})} = \underset{P_{1}, P_{2} \in
		\mathcal{P}}{\max} |P_{1}P_{2}|
		.
	\end{equation}
\end{dfn}

Every integral point set has also a characteristic~\cite{kemnitz1988punktmengen,kurz2005characteristic}.

\begin{dfn}
	The characteristic of a PIPS $\mathcal{P}$ is a squarefree number $q$
	such that the area of any triangle $M_1M_2M_3$, $\{M_1,M_2,M_3\}\subset \mathcal{P}$,
	is commensurable with $\sqrt{q}$.
\end{dfn}

The following notion was introduced in~\cite{kurz2008bounds}.

\begin{dfn}
	The function $d(m, n)$ is the minimal possible diameter of
	an integral point set $\mathcal{P}$ of $n$ points in
	$m$-dimensional Euclidean space $\mathbb{R}^{m}$.
\end{dfn}

The computation of exact values of $d(m,n)$ is a difficult problem.
In 2003, Solymosi proved~\cite{solymosi2003note} that $d(2,n) \geq cn$ for some constant $c$;
now it is known~\cite{my-pps-linear-bound-2019} that $c>5/11$.
The most recent advances for higher dimensions can be found in~\cite{nozaki2013lower}.
So, even finding the estimates is not easy.

A long list of known bounds for $d(m, n)$ and some exact values
can be found in~\cite[Theorem 1]{kurz2008bounds} or in~\cite{our-vmmsh-2018}.
Below we will discuss the following ones presented at~\cite{kurz2008bounds}:
\begin{align}
	d(m, 2m + 1) \leq 8  \label{eq:d_m_2m+1}\\
	d(m, 2m + 2) \leq 13 \label{eq:d_m_2m+2}\\
	d(m, 3m) \leq 109    \label{eq:d_m_3m}
\end{align}
and the following theorem \cite[Theorem 2.1]{kurz2008bounds}.

\begin{theorem}
	\label{thm:Kurz_blowup}
	Let $\mathcal{P}$ be a planar integral point set consisting of
	$n - 2$ points on the line $l_1$ and two points $P_{1}$ and $P_{2}$ on a
	parallel line $l_2$ with distance $r$ between $l_{1}$ and $l_{2}$. If there
	exist positive integers $v$, $w$ with $f^{2} + v^{2}
	= w^{2}$ and $v < 2r$, where $|P_{1}P_{2}| = f$,
	then

	\begin{equation}\label{formula1}
		d(m, n + 2(m - 2)) \leq \max(w, \operatorname{diam}(\mathcal{P}))
	\end{equation}

\end{theorem}

In this paper we further study the various features of IPS.
In Section~\ref{sec:classif}, we discuss the classification of planar integral points sets;
in Sections~\ref{sec:bounds_rails}~and~\ref{sec:bounds_pyramid}, we present some upper bounds for $d(m,n)$
based on planar integral point sets of special types
and provide some general constructions of such bounds.
All the bounds, like~\eqref{eq:d_m_2m+1}--\eqref{eq:d_m_3m}, are of the form $d(m, km+p)\leq c_{k,p}$.

\section{Classification of planar integral point sets}
\label{sec:classif}

\subsection{Integral point sets situated on two straight lines}

\begin{dfn}
	A planar integral point sets of $n$ points with $n-1$ points on a straight line is called
	a \textit{facher} set.
\end{dfn}
Most of the examples of planar integral point sets turns to be facher sets.
In~\cite{antonov2008maximal}, facher sets with characteristic 1 are called \textit{semi-crabs}.
For $9 \leq n \leq 122$, the diameter $d(2,n)$ is attained on a facher point set~\cite{kurz2008minimum}.

Non-facher integral point sets situated on two straight lines
can be further classificated in the following three cases:

\begin{dfn}
	A planar integral point set situated on two parallel straight lines
	is called a \textit{rails} set.
\end{dfn}

Among the rails sets, sets with 2 points on one line and all the other on another line prevail.
In Section~\ref{sec:bounds_rails} below, some examples of such sets are discussed;
they include an example of rails sets with 2 points on one line and 38 point on the other line
(Fig.~\ref{picture_40.png})
and
an example of rails set with 4 points on one line and 4 points on the other
(Fig.~\ref{picture.png}).
We still do not know whether there are rails set with 4 points on one line and 5 points on the other;
and the same for rails sets with 3 points on one line and 9 points on the other.

For the sake of brevity, we will hereafter use the notation
$\sqrt{p}/q * \{ (x_1,y_1), ...,$ $ (x_n, y_n)  \}$
from~\cite{our-ped-2018,our-pmm-2018,our-vmmsh-2018},
which means that each abscissa is multiplied by $1/q$
and each ordinate is multiplied by $\sqrt{p}/q$,  that is
$$
	\sqrt{p}/q * \{ (x_1,y_1), ..., (x_n, y_n)  \}
	=
	\left\{ \left(\frac{x_1}{q},\frac{y_1\sqrt{p}}{q}\right), ..., \left(\frac{x_n}{q},   \frac{y_n\sqrt{p}}{q}\right)  \right\}.
$$

\begin{figure}[h!]
	\center{\includegraphics[width=0.9\linewidth]{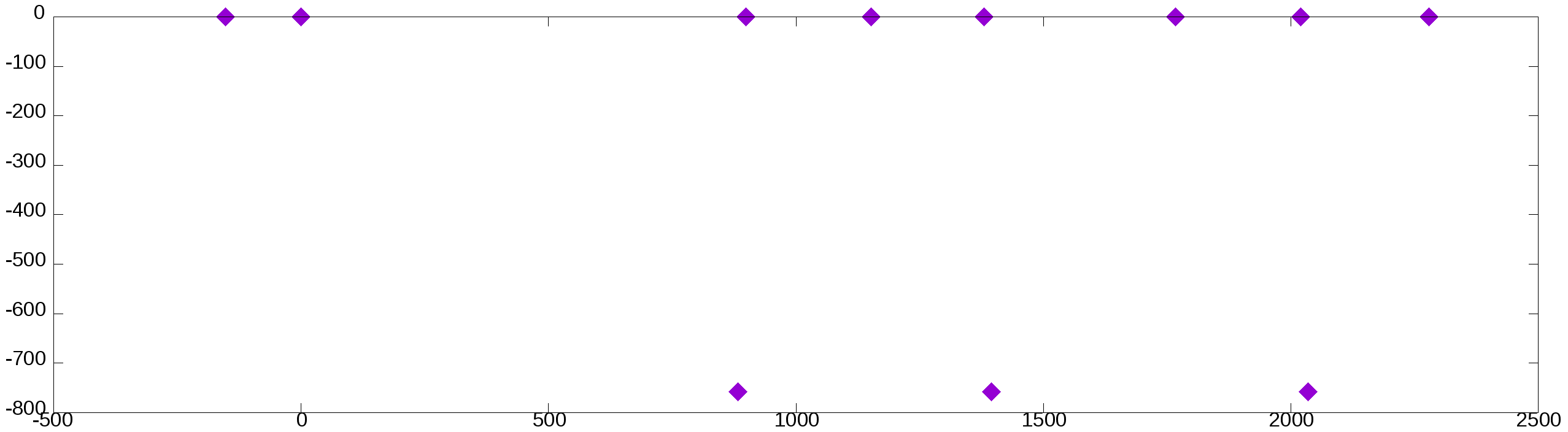}}
	\parbox{0.9\linewidth}{\caption{PIPS of cardinality 11 and diameter 2423}
	\label{11_2432_255255_a0a1b7f805cfc93080b4608505a0f45d_rails_3.png}}
\end{figure}

On Fig.~\ref{11_2432_255255_a0a1b7f805cfc93080b4608505a0f45d_rails_3.png} we show the following rails set with 3 points on one line and 8 points on the other:
\begin{multline*}
\mathcal{P}_{3,8}=\sqrt{255255}/2*\{
	( 1767 ; -3);
	( 2791 ; -3);
	( 4071 ; -3);
\\
	( -306 ; 0);
	( 0 ; 0);
	( 1798 ; 0);
	( 2304 ; 0);
	( 2760 ; 0);
	( 3534 ; 0);
	( 4040 ; 0);
	( 4558 ; 0)
\}
\end{multline*}
It's notable that $\operatorname{char} \mathcal{P}_{3,8}=255255 = 3 \cdot 5 \cdot 7 \cdot 11 \cdot 13 \cdot 17$
is a product of the first 6 odd primes,
and the distances between the three points on the lower line are $512=2^{9}$, $640=2^7 \cdot 5$ and
$512+640 = 1152 = 2^7 \cdot 3^2$.

\begin{figure}[h!]
	\center{\includegraphics[width=0.9\linewidth]{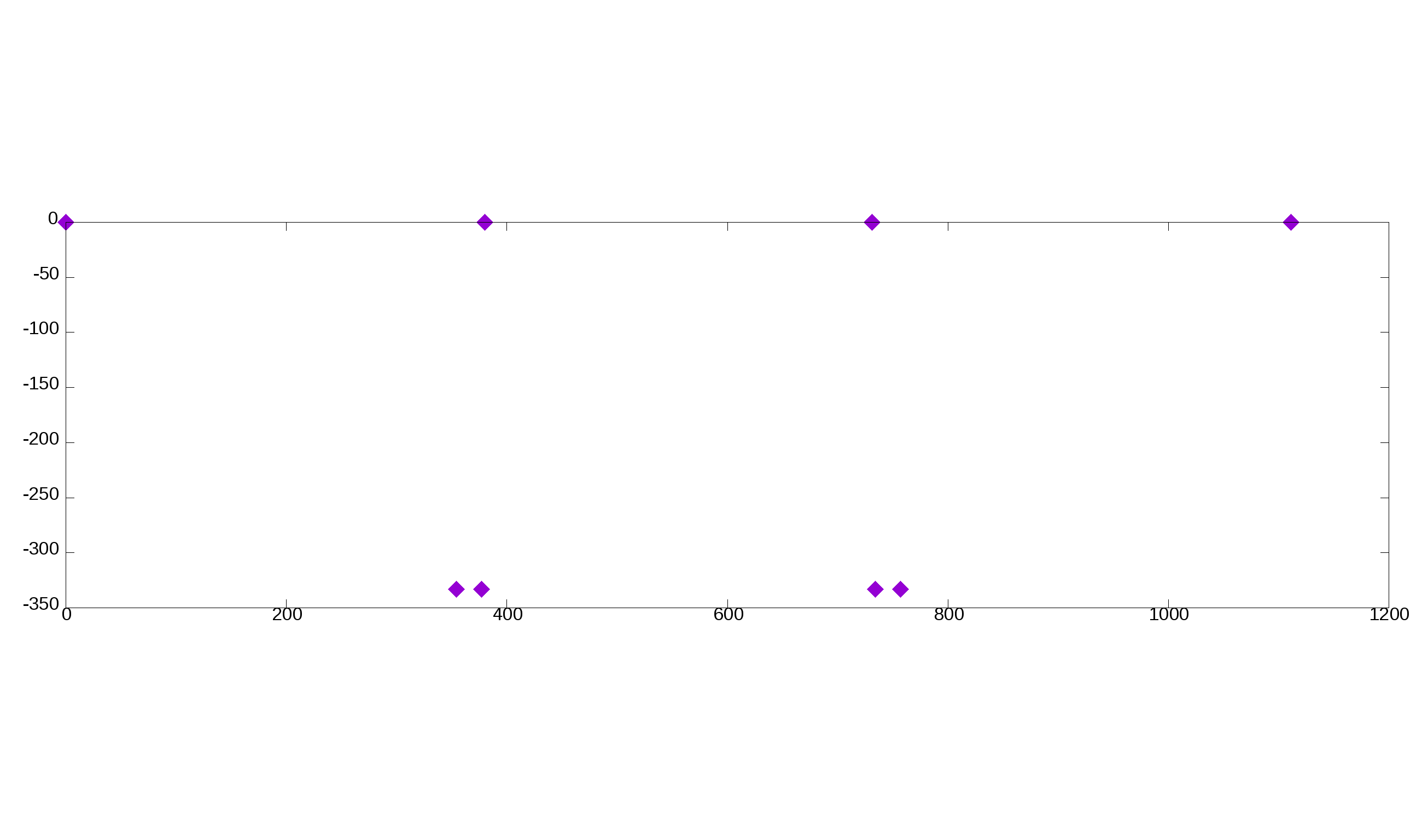}}
	\parbox{0.9\linewidth}{\caption{PIPS of cardinality 8 and diameter 1111}
	\label{8_1111_770_5666c83a03a80e1d3f8f8e8186a145f6__rails_4.png}}
\end{figure}

Rails sets with 4 and 4 points on the lines are somewhat rather common;
an example of such a set is shown on Fig.~\ref{8_1111_770_5666c83a03a80e1d3f8f8e8186a145f6__rails_4.png}:
\begin{equation*}
	\sqrt{770}/1*\{
		( 0 ; 0);
		( 380 ; 0);
		( 731 ; 0);
		( 1111 ; 0);
		( 354 ; -12);
		( 377 ; -12);
		( 734 ; -12);
		( 757 ; -12)
	\}
\end{equation*}

\begin{dfn}
	A planar integral point set situated on two perpendicular straight lines
	is called a \textit{cross} set.
\end{dfn}
Every cross set has characteristic 1;
in~\cite{antonov2008maximal}, cross sets with only 2 points out of one of the lines that have two axes of symmetry are called \textit{crabs}.

\begin{dfn}
	A planar integral point set situated on two straight lines
	that are not parallel nor perpendicular,
	is called a \textit{scissors} set.
\end{dfn}

There is an important subclass of scissors sets.

\begin{dfn}
	A scissors set with an axis of symmetry,
	which is the angle bisector for the straight lines,
	is called a \textit{pyramid} set.
\end{dfn}

Examples of pyramid sets are given in Section~\ref{sec:bounds_pyramid} on
Fig.~\ref{picture_56.png} and Fig.~\ref{picture_1260.png}.
We don't have any examples of pyramid sets with cardinality more than 9.

It still remains unknown whether a scissors set can have another axis of symmetry;
on Fig.~\ref{15_3364_91_30e991f82a54fb9caeca5b6a6776b981_X.png} the following scissors set with central symmetry is shown:
\begin{multline*}
	\sqrt{91}/1*\{
	( 0 ; 0);
	( 792 ; 120);
	( -792 ; -120);
\\
	( \pm256 ; 0);
	( \pm646 ; 0);
	( \pm693 ; 0);
	( \pm1152 ; 0);
	( \pm1242 ; 0);
	( \pm1682 ; 0)
\}
\end{multline*}

\subsection{Other integral point sets}

\begin{dfn}
	A planar integral point set that is situated on a circle is called a \textit{circular}
	point set.
\end{dfn}

Circular sets are very important examples
of integral point sets~\cite{harborth1993upper,piepmeyer1996maximum,bat2018number};
for example, the best known upper bound for $d(2,n)$
is attained on a circular set.

\begin{dfn}
	A planar integral point set that is situated on the conjunction of a circle with its center,
	is called a \textit{centered-circular} point set.
\end{dfn}

\begin{remark}
	Every centered-circular point set has characteristic 1.
	Indeed, let $\mathcal{P} = \{M_0, M_1,$ $ M_2, M_3, ..., M_k\}$ be a centered-circular point set
	with $M_0$ being the center.
	by the Law of sines, in the triangle $M_1 M_2 M_3$
	\begin{equation}
		\sin \angle M_1 M_2 M_3 = \frac{|M_1M_3|}{2R}
		,
	\end{equation}
	where $R = |M_0 M_1| $ is the radius of the triangle's circumcircle.
	Then the area of the triangle $M_1 M_2 M_3$ is
	\begin{equation}
		S = \frac12 |M_2 M_1| \cdot |M_2 M_3| \cdot \sin \angle M_1 M_2 M_3 =
		\frac12 |M_2 M_1| \cdot |M_2 M_3| \cdot \frac{|M_1M_3|}{2|M_0 M_1|} \in\mathbb Q
		.
	\end{equation}
\end{remark}

\begin{figure}[htbp]
	\includegraphics[width=.42\linewidth]{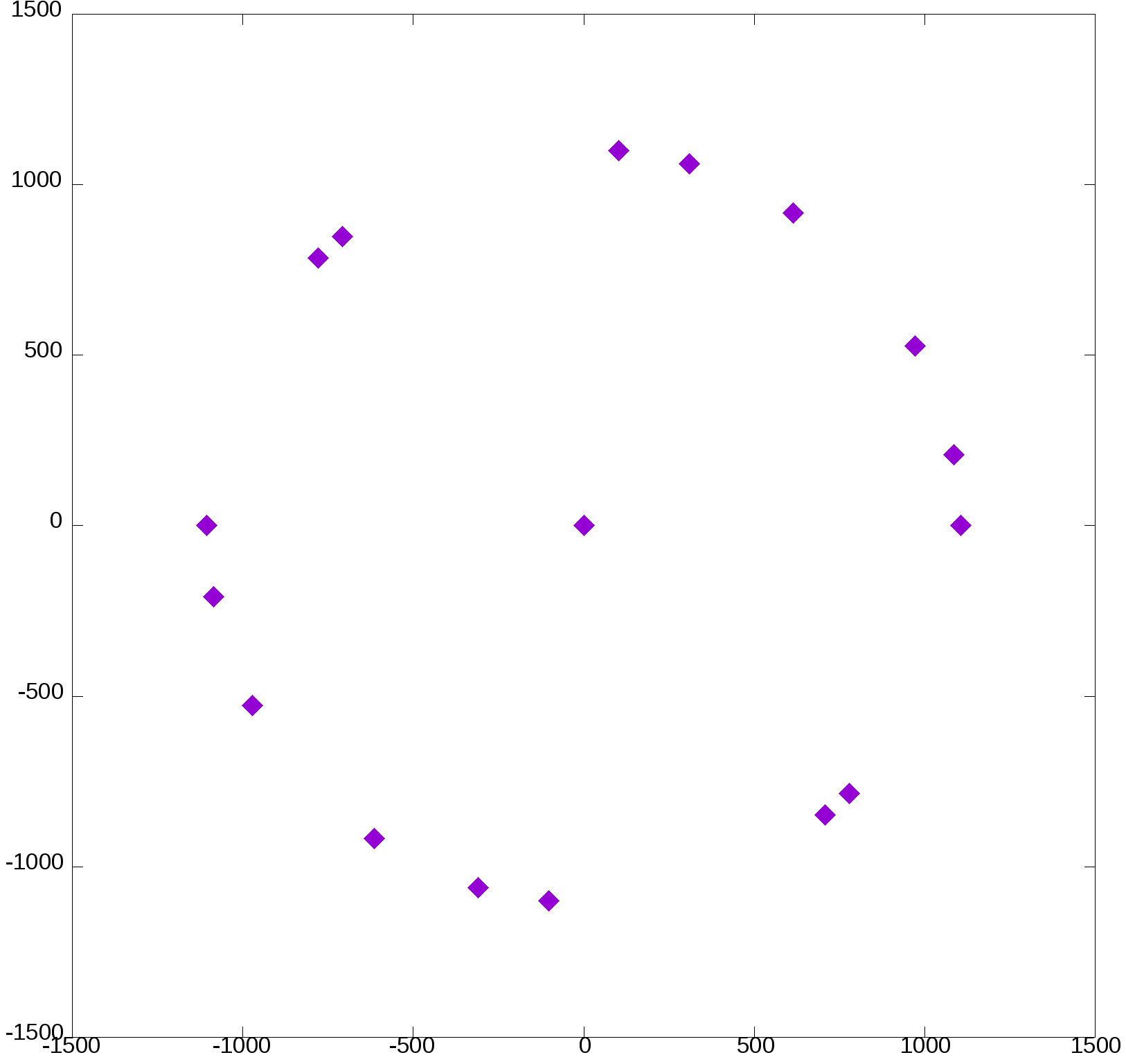}
	\hfill
	\includegraphics[width=.545\linewidth]{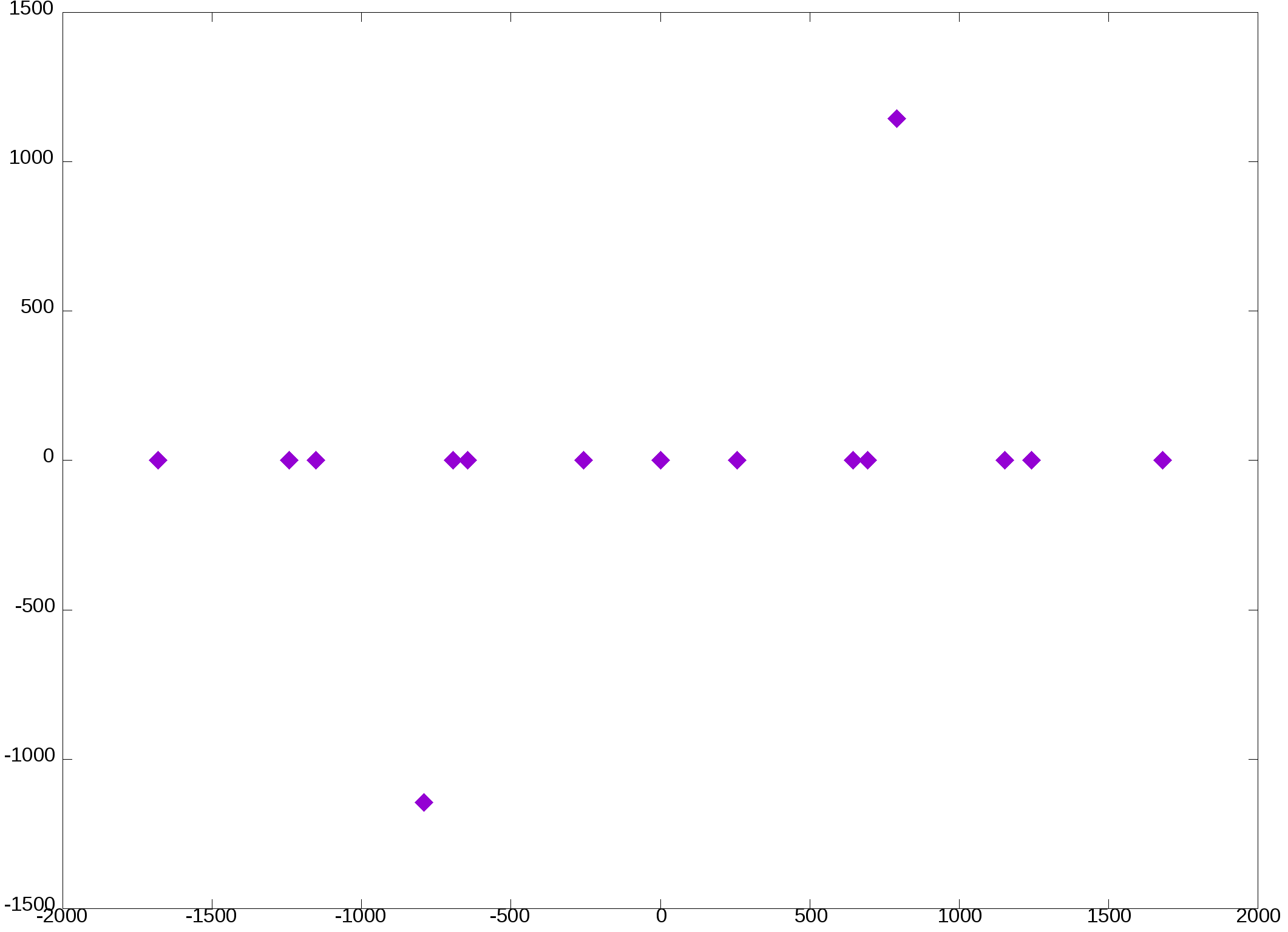}
	\\
	\parbox{.37\linewidth}{\caption{A centered-circular PIPS of cardinality 17 and diameter 2210}
	\label{17_2210_1_acc4a56887abb3c3a315712ea083f0e7_centered_circle.png}}
	\hfill
	\parbox{.5\linewidth}{\caption{A scissors PIPS of cardinality 15 and diameter 3364}
	\label{15_3364_91_30e991f82a54fb9caeca5b6a6776b981_X.png}}
\end{figure}

On Fig.~\ref{17_2210_1_acc4a56887abb3c3a315712ea083f0e7_centered_circle.png} the following centered-circular PIPS is shown:
\begin{multline*}
	\sqrt{1}/1105*\{
	( 0 ; 0);
	( \pm1221025 ; 0);
	( 113953 ; 1215696);
	( -113953 ; -1215696);
	\\
	( 341887 ; 1172184);
	( -341887 ; -1172184);
	( 680225 ; 1014000);
	( -680225 ; -1014000);
	\\
	( 782977 ; -936936);
	( -782977 ; 936936);
	( 859775 ; -867000);
	( -859775 ; 867000);
	\\
	( 1073057 ; 582624);
	( -1073057 ; -582624);
	( 1198975 ; 231000);
	( -1198975 ; -231000);
\}
\end{multline*}

Most of the known planar integral point sets fall into these six classes;
however, there are some sophisticated examples which do not~\cite{kreisel2008heptagon,kurz2013constructing,avdeev2021particular}.
A lower bound for the diameter of a PIPS with no collinear triples (so-called PIPS \textit{in semi-general position})
is presented in~\cite{my-semi-general-5-4-bound-2019}.


\section{Bounds based on rails sets}
\label{sec:bounds_rails}

Let $i = \overline{1, k}$ denote the enumeration of all $i$
from $1$ to $k$.
Theorem \cite[Theorem 2.1]{kurz2008bounds} can be generalized as follows:

\begin{theorem}
	\label{thm:rails_blowup}
	Let $\mathcal{P}$ be a planar integral point set consisting of
	$n - k$ points on the line $l_1$ and $k$ points $P_1$, $P_2$, ..., $P_k$ on a
	parallel line $l_2$ with distance $r$ between $l_1$ and $l_2$. If there
	exist positive integers $v$, $w_{ij}$, with $f_{ij}^{2} + v^{2}
	= w_{ij}^{2}$ and $v < 2r$, where $i = \overline{1, k - 1}$, $j =
	\overline{i + 1, k}$, $|P_{i}P_{j}| = f_{ij}$,
	then

	\begin{equation}
		d(m, n + k(m - 2)) \leq \max(w_{1k}, \operatorname{diam(\mathcal{P})})
	\end{equation}

\end{theorem}

Below we give the examples of planar integral point sets and the corresponding
estimates of the function $d(m, n)$ for $n = 2m + k$, $3 \leq k \leq 36$.

\begin{figure}[htbp]
	\includegraphics[width=.48\linewidth]{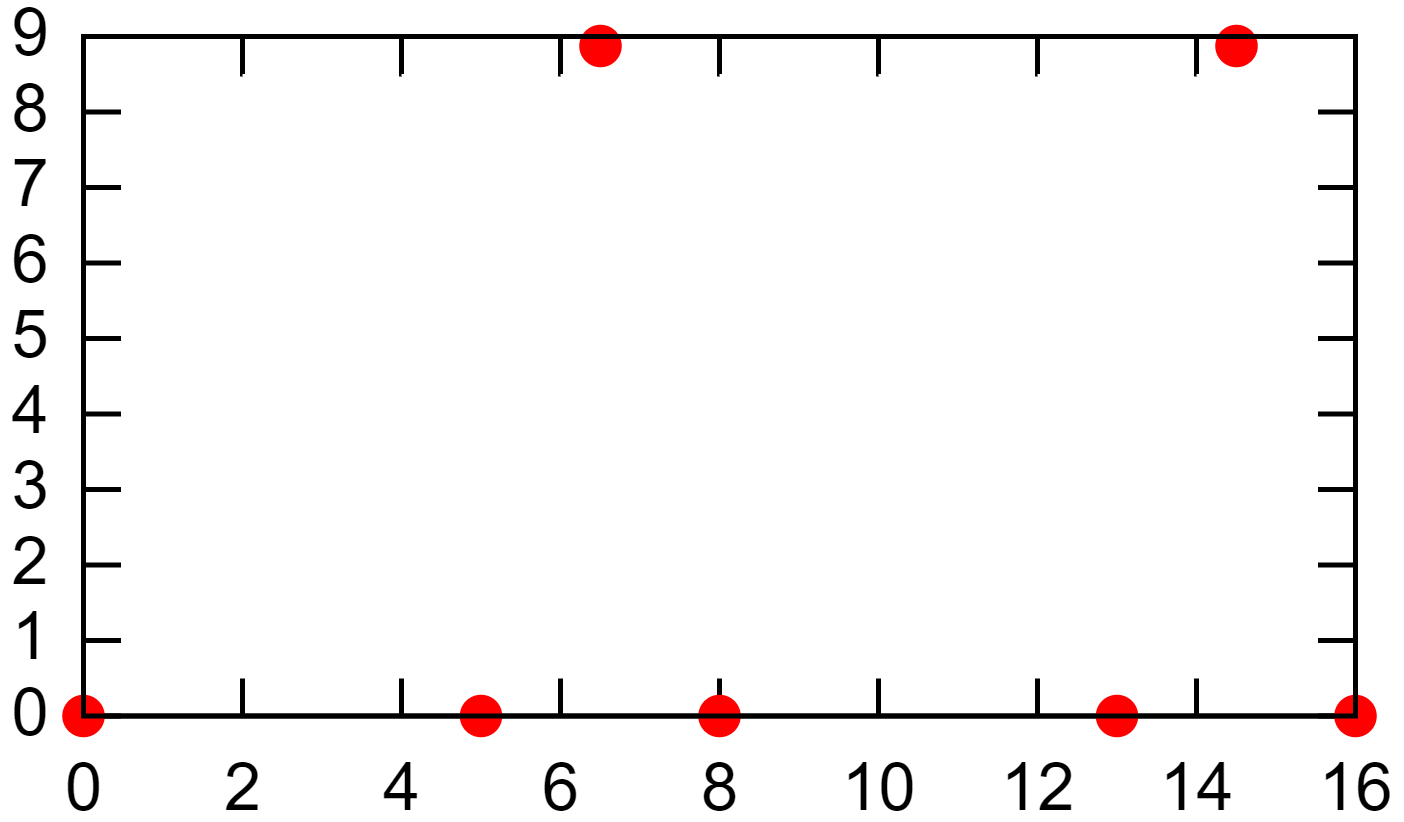}
	\hfill
	\includegraphics[width=.48\linewidth]{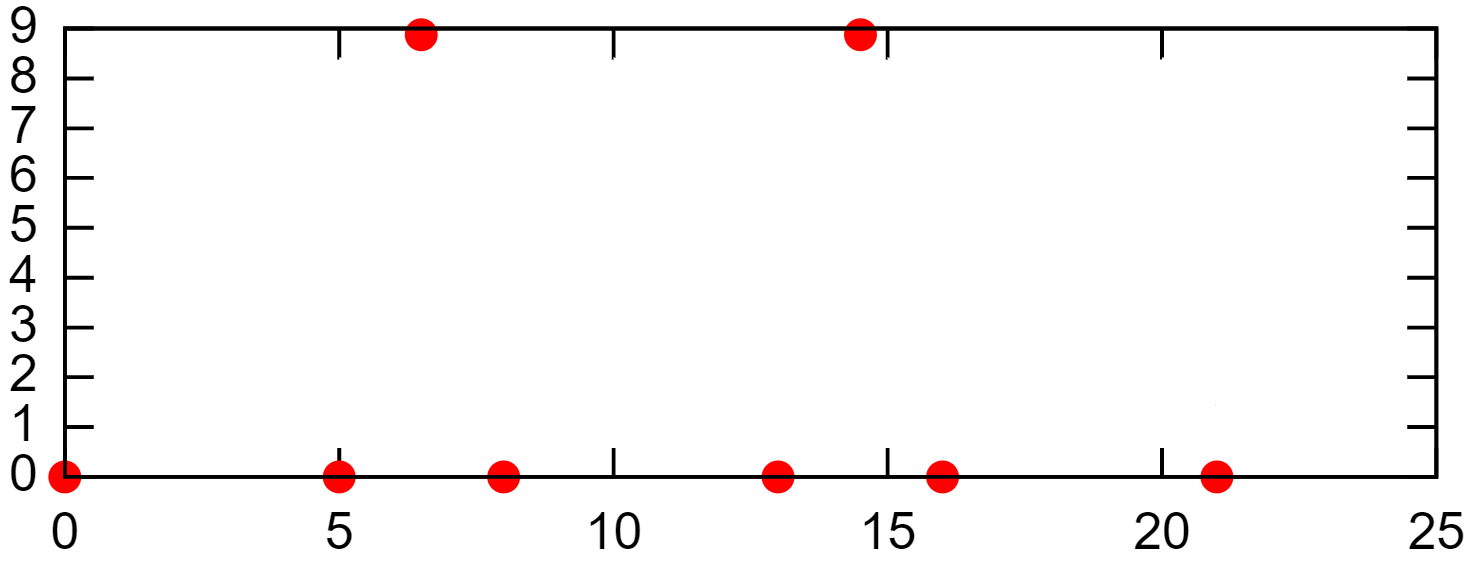}
	\\
	\parbox{.48\linewidth}{\caption{PIPS of cardinality 7 and diameter 17}
	\label{picture_d17.png}}
	\hfill
	\parbox{.48\linewidth}{\caption{PIPS of cardinality 8 and diameter 21}
	\label{picture_d21.png}}
\end{figure}

Below $f$ and $w$ are those from Theorem~\ref{thm:Kurz_blowup}.

\begin{itemize}
\setlength{\itemsep}{-1mm}

\item
$\mathcal{P}=\sqrt{315}/{2} * \{ (13, 1),
(29, 1),
(0, 0),
(10, 0),
(16, 0),
(26, 0),
(32, 0)\}
$
(Figure~\ref{picture_d17.png})

$f = 8$, $v = 6$, $w = 10$, $\operatorname{diam(\mathcal{P})} = 17$,

which gives $d(m, 2m + 3) \leq 17$.

\item
$\mathcal{P}=\sqrt{315}/{2} * \{ (13, 1),
(29, 1),
(0, 0),
(10, 0),
(16, 0),
(26, 0),
(32, 0),
(42, 0)\}
$
(Figure~\ref{picture_d21.png})

$f = 8$, $v = 6$, $w = 10$, $\operatorname{diam(\mathcal{P})} = 21$,

which gives $d(m, 2m + 4) \leq 21$.

\item
$\mathcal{P}=\sqrt{70}/{1} * \{ (\pm 44, 12),
(-78 , 0),
(\pm 62, 0),
(\pm 55 , 0),
(\pm 10 , 0)\}
$
(Figure~\ref{picture_2.png})

$f = 88$, $v = 66$, $w = 110$, $\operatorname{diam(\mathcal{P})} = 158$,

which gives $d(m, 2m + 5) \leq 158$.

\item
$\mathcal{P}=\sqrt{70}/{1} * \{ (\pm 44, 12),
(\pm 78 , 0),
(\pm 62, 0),
(\pm 55 , 0),
(\pm 10 , 0)\}
$
(Figure~\ref{picture_3.png}).

$f = 88$, $v = 66$, $w = 110$, $\operatorname{diam(\mathcal{P})} = 158$,

which gives $d(m, 2m + 6) \leq 158$.

\begin{figure}[h!]
	\includegraphics[width=.48\linewidth]{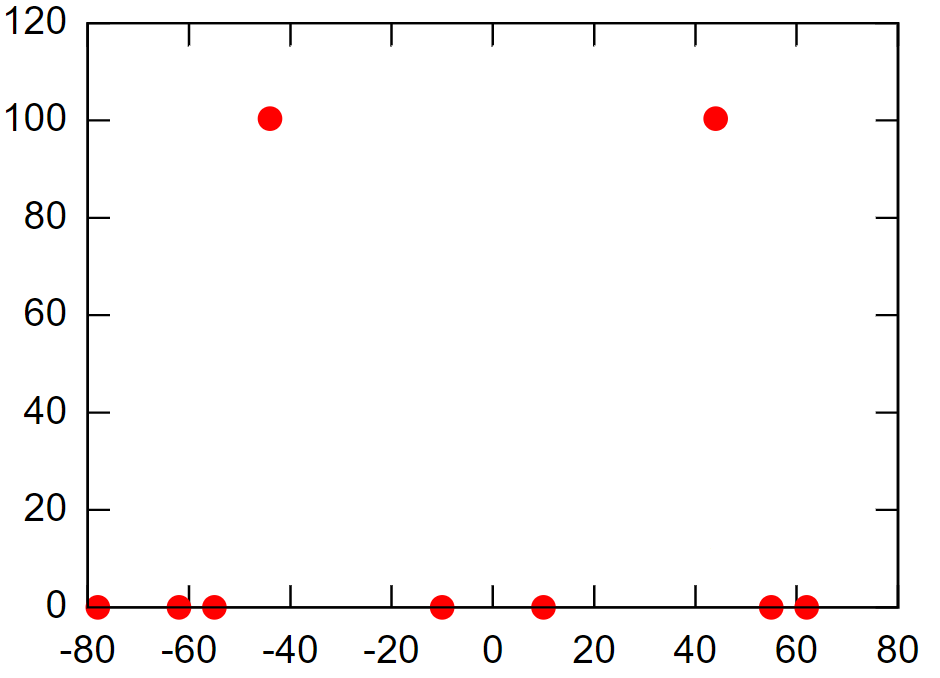}
	\hfill
	\includegraphics[width=.48\linewidth]{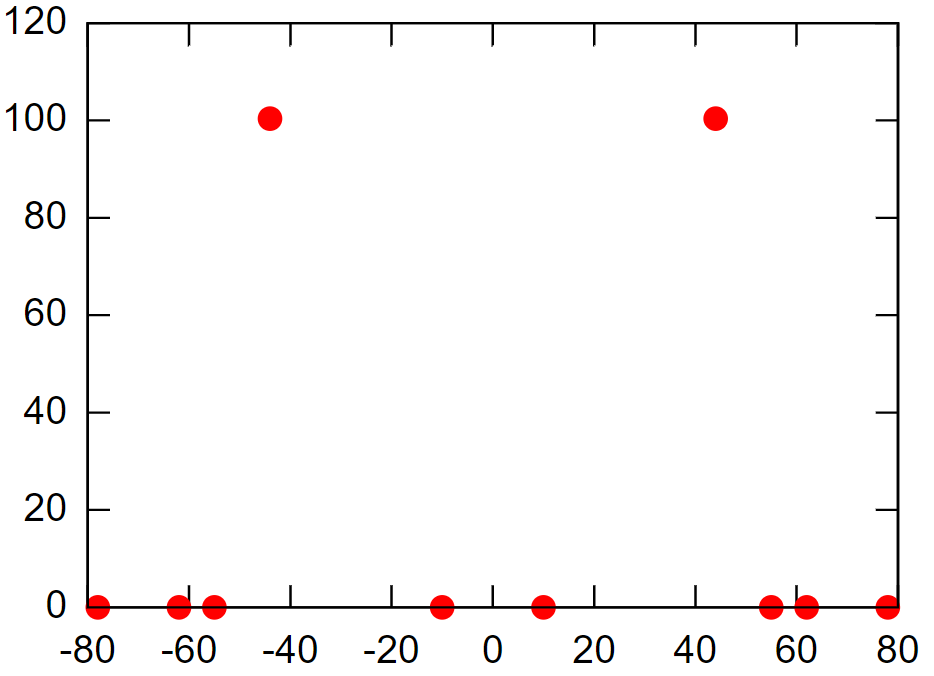}
	\\
	\parbox{.48\linewidth}{\caption{PIPS of cardinality 9 and diameter 158}
	\label{picture_2.png}}
	\hfill
	\parbox{.48\linewidth}{\caption{PIPS of cardinality 10 and diameter 158}
	\label{picture_3.png}}
\end{figure}

\item
$\mathcal{P}=\sqrt{70}/{2} * \{ (\pm 99, 24),
(\pm 207 , 0),
(\pm 145 , 0),
(\pm 63 , 0),
(\pm 25 , 0),
(-297 , 0)\}
$
(Figure~\ref{picture_4.png}).

$f = 99$, $v = 20$, $w = 101$, $\operatorname{diam(\mathcal{P})} = 252$,

which gives $d(m, 2m + 7) \leq 252$.

\begin{figure}[htbp]
	\includegraphics[width=.48\linewidth]{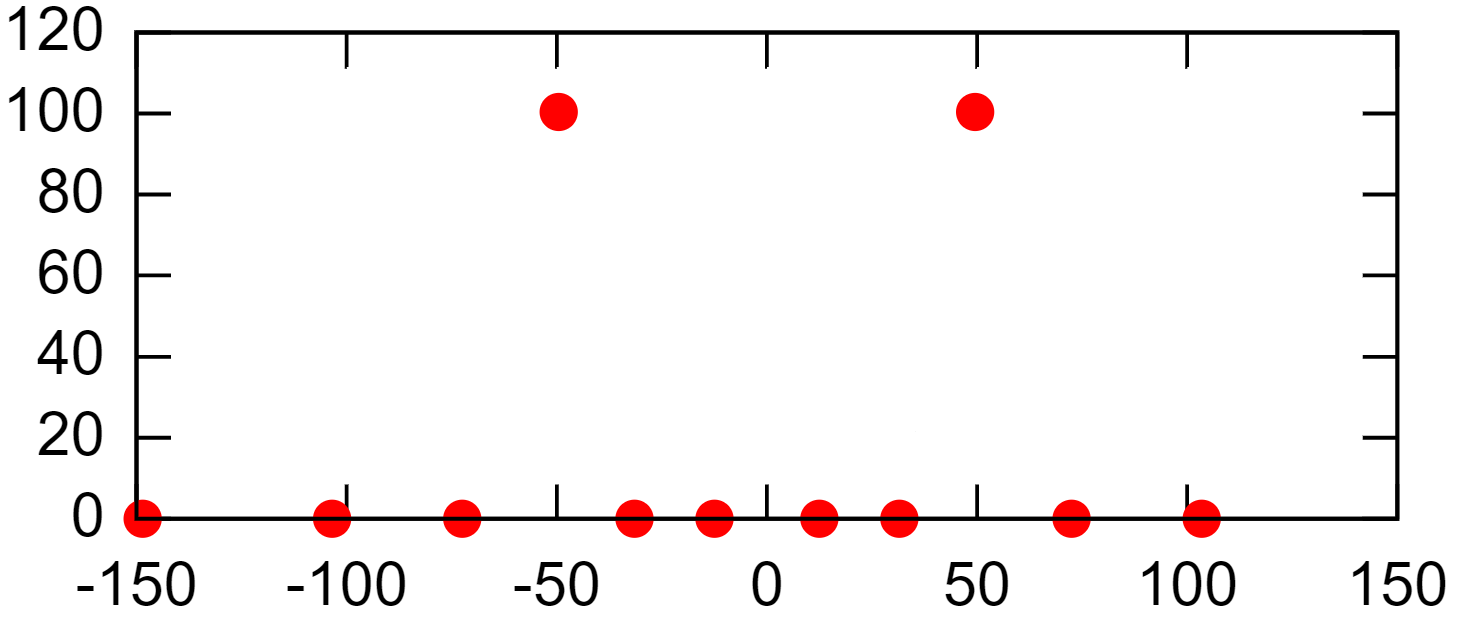}
	\hfill
	\includegraphics[width=.48\linewidth]{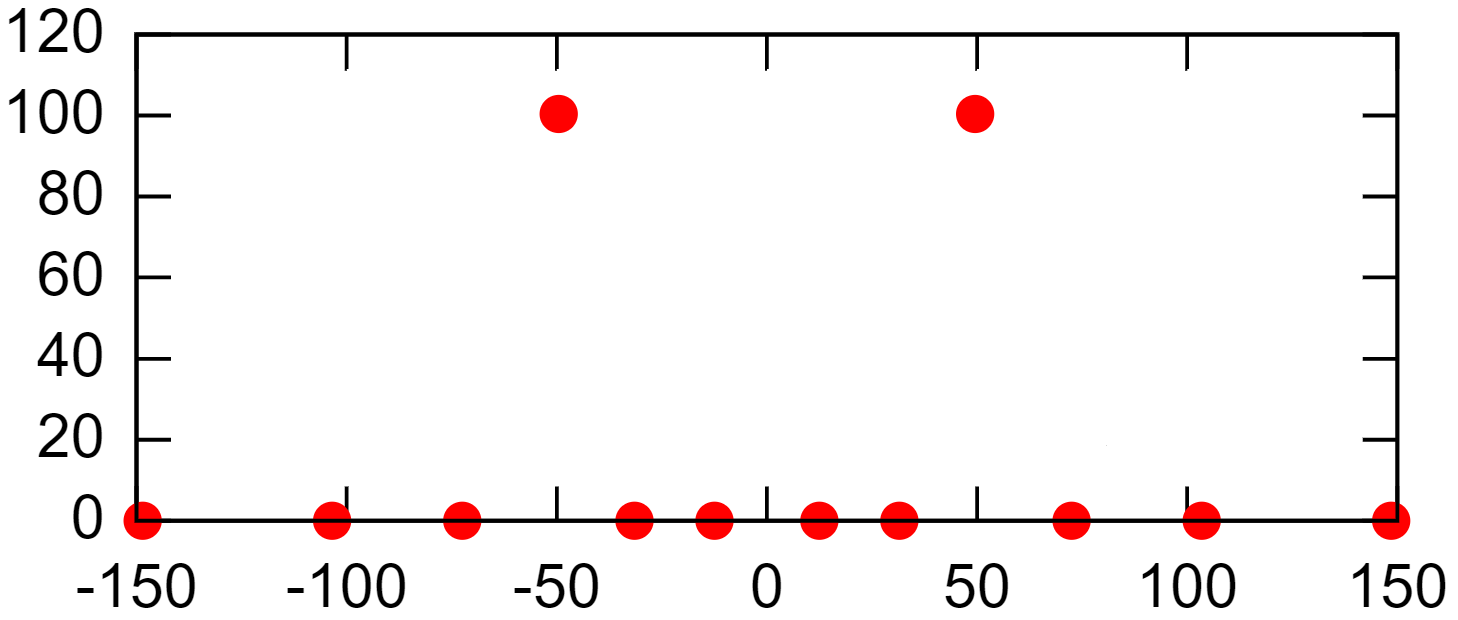}
	\\
	\parbox{.48\linewidth}{\caption{PIPS of cardinality 11 and diameter 252}
	\label{picture_4.png}}
	\hfill
	\parbox{.48\linewidth}{\caption{PIPS of cardinality 12 and diameter 297}
	\label{picture_5.png}}
\end{figure}

\item
$\mathcal{P}=\sqrt{70}/{2} * \{ (\pm 99, 24),
(\pm 297 , 0),
(\pm 207 , 0),
(\pm 145 , 0),
(\pm 63 , 0),
(\pm 25 , 0)\}
$
(Figure~\ref{picture_5.png}).

$f = 99$, $v = 20$, $w = 101$, $\operatorname{diam(\mathcal{P})} = 297$,

which gives $d(m, 2m + 8) \leq 297$.

\item
$\mathcal{P}=\sqrt{19019}/{2} * \{ (\pm 200, 3),
(\pm 873 , 0),
(\pm 615 , 0),
(\pm 377 , 0),
(\pm 215 , 0),
(\pm 23 , 0),
$

$
(-1273 , 0)\}
$
(Figure~\ref{picture_6.png}).

$f = 200$, $v = 45$, $w = 205$, $\operatorname{diam(\mathcal{P})} = 1073$,

which gives $d(m, 2m + 9) \leq 1073$.

\begin{figure}[h!]
\center{\includegraphics[width=0.7\linewidth]{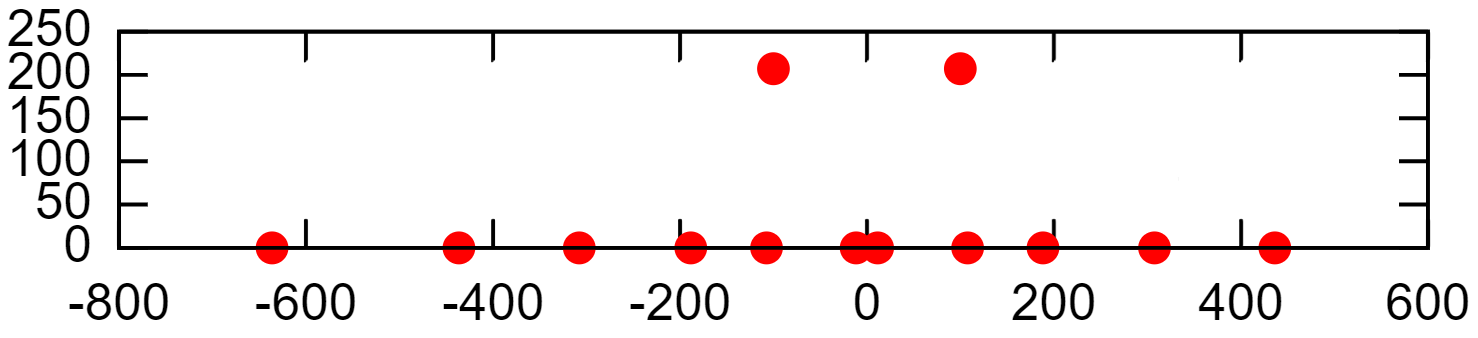}}
\parbox{0.7\linewidth}{\caption{PIPS of cardinality 13 and diameter 1073}
\label{picture_6.png}}
\end{figure}

\item
$\mathcal{P}=\sqrt{19019}/{2} * \{ (\pm 200, 3),
(\pm 1273 , 0),
(\pm 873 , 0),
(\pm 615 , 0),
(\pm 377 , 0),
(\pm 215 , 0),
$

$
(\pm 23 , 0)\}
$
(Figure~\ref{picture_7.png}).

$f = 200$, $v = 45$, $w = 205$, $\operatorname{diam(\mathcal{P})} = 1273$,

which gives $d(m, 2m + 10) \leq 1273$.

\begin{figure}[h!]
\center{\includegraphics[width=0.7\linewidth]{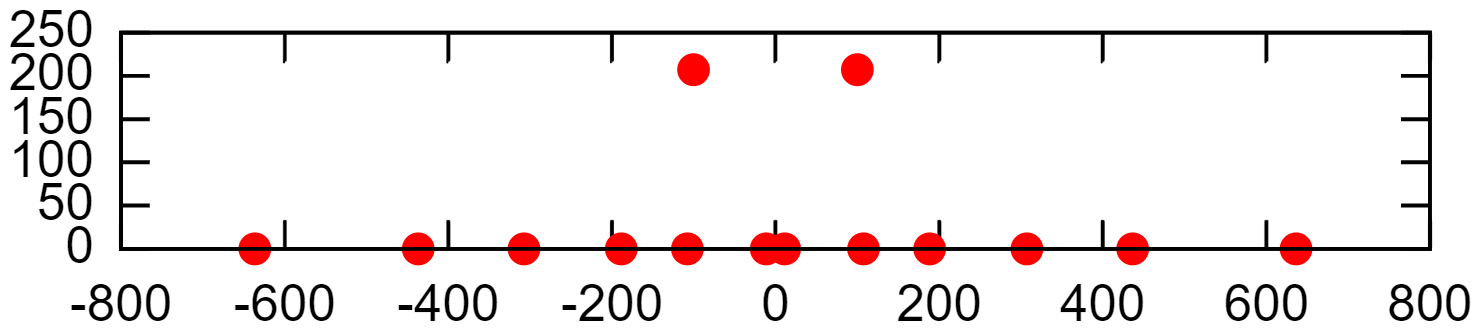}}
\parbox{0.7\linewidth}{\caption{PIPS of cardinality 14 and diameter 1273}
\label{picture_7.png}}
\end{figure}

\item
$\mathcal{P}=\sqrt{385}/{2} * \{ (\pm 1105, 48),
(\pm 1587 , 0),
(\pm 1269 , 0),
(\pm 763 , 0),
(\pm 623 , 0),
(\pm 529 , 0),
$

$
(\pm 339 , 0),
(-2189 , 0)\}
$
(Figure~\ref{picture_8.png}).

$f = 1105$, $v = 300$, $w = 1145$, $\operatorname{diam(\mathcal{P})} = 1888$,

which gives $d(m, 2m + 11) \leq 1888$.

\begin{figure}[h!]
\center{\includegraphics[width=0.75\linewidth]{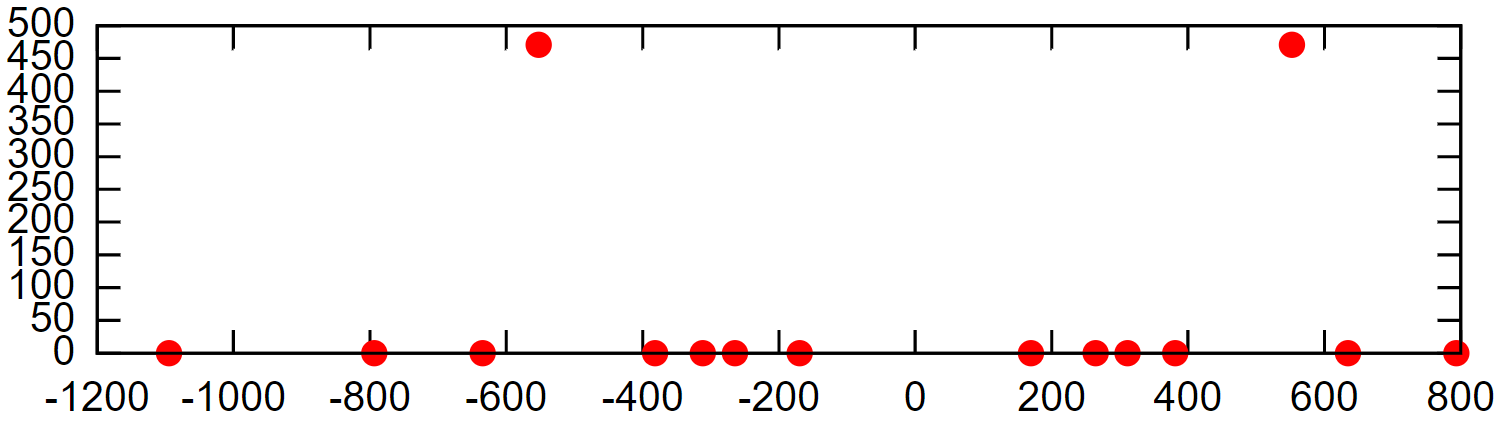}}
\parbox{0.75\linewidth}{\caption{PIPS of cardinality 15 and diameter 1888}
\label{picture_8.png}}
\end{figure}

\item
$\mathcal{P}=\sqrt{385}/{2} * \{ (\pm 1105, 48),
(\pm 2189 , 0),
(\pm 1587 , 0),
(\pm 1269 , 0),
(\pm 763 , 0),
(\pm 623 , 0),
$

$
(\pm 529 , 0),
(\pm 339 , 0)\}
$
(Figure~\ref{picture_9.png}).

$f = 1105$, $v = 300$, $w = 1145$, $\operatorname{diam(\mathcal{P})} = 2189$,

which gives $d(m, 2m + 12) \leq 2189$.

\begin{figure}[h!]
\center{\includegraphics[width=0.75\linewidth]{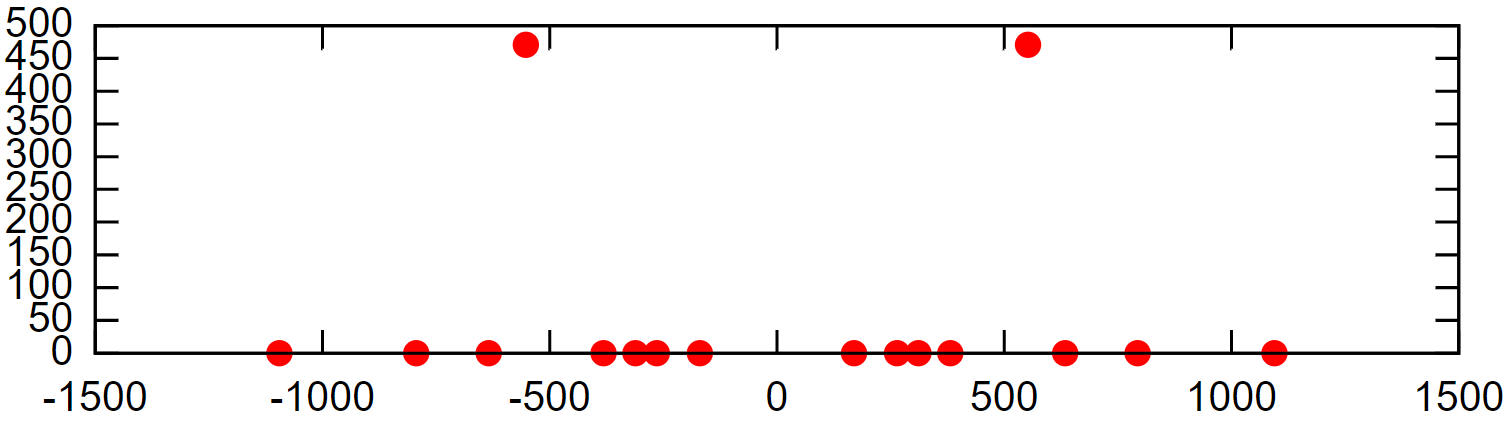}}
\parbox{0.75\linewidth}{\caption{PIPS of cardinality 16 and diameter 2189}
\label{picture_9.png}}
\end{figure}

\item
$\mathcal{P}=\sqrt{154}/{1} * \{ (\pm 874, 60),
(\pm 1376 , 0),
(\pm 1036 , 0),
(\pm 899 , 0),
(\pm 849 , 0),
(\pm 613 , 0),
$

$
(\pm 576 , 0),
(\pm 100 , 0),
(-1848 , 0)\}
$
(Figure~\ref{picture_13.png}).

$f = 1748$, $v = 336$, $w = 1780$, $\operatorname{diam(\mathcal{P})} = 3224$,

which gives $d(m, 2m + 13) \leq 3224$.

\begin{figure}[h!]
\center{\includegraphics[width=0.85\linewidth]{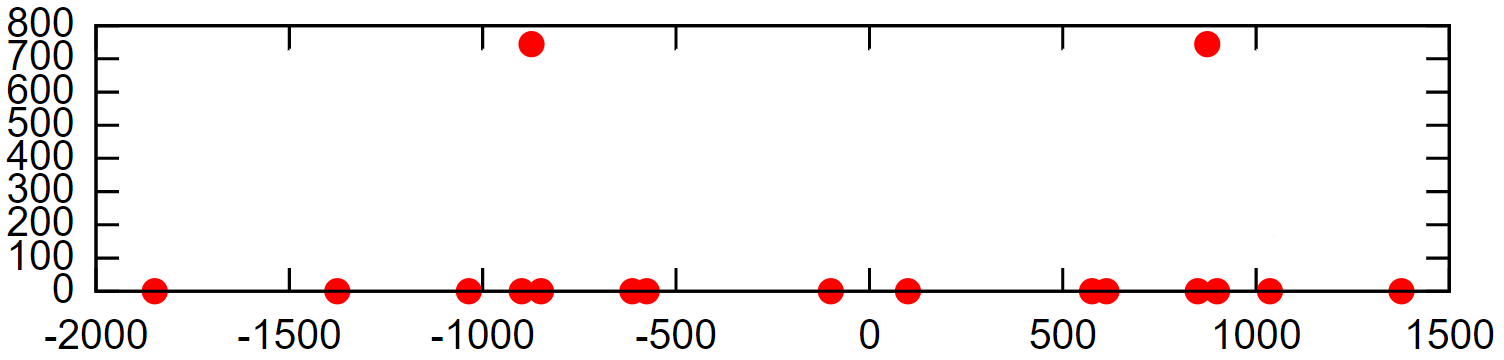}}
\parbox{0.85\linewidth}{\caption{PIPS of cardinality 17 and diameter 3224}
\label{picture_13.png}}
\end{figure}

\item
$\mathcal{P}=\sqrt{154}/{1} * \{ (\pm 874, 60),
(\pm 1848 , 0),
(\pm 1376 , 0),
(\pm 1036 , 0),
(\pm 899 , 0),
(\pm 849 , 0),
$

$
(\pm 613 , 0),
(\pm 576 , 0),
(\pm 100 , 0)\}
$
(Figure~\ref{picture_14.png}).

$f = 1748$, $v = 336$, $w = 1780$, $\operatorname{diam(\mathcal{P})} = 3696$,

which gives $d(m, 2m + 14) \leq 3696$.

\begin{figure}[h!]
\center{\includegraphics[width=0.85\linewidth]{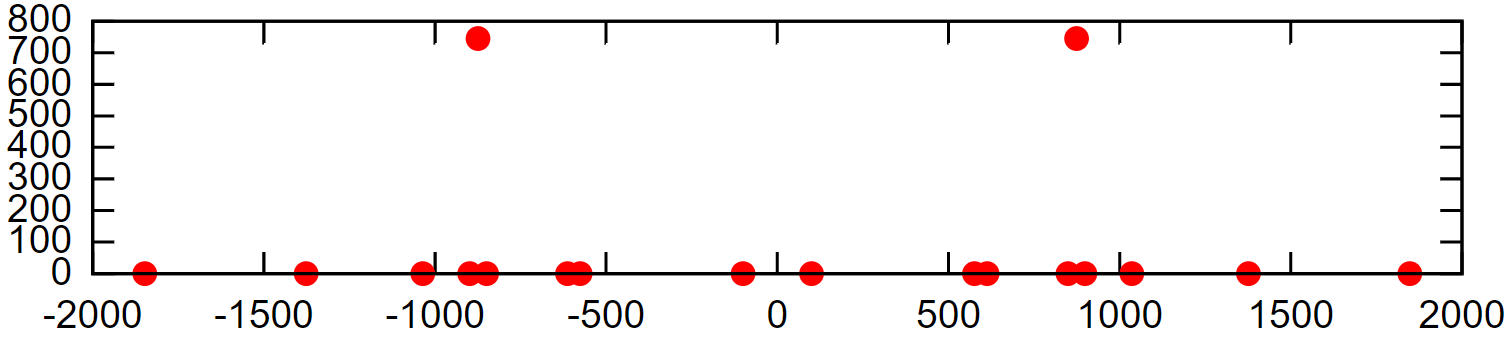}}
\parbox{0.85\linewidth}{\caption{PIPS of cardinality 18 and diameter 3696}
\label{picture_14.png}}
\end{figure}

\item
$\mathcal{P}=\sqrt{154}/{1} * \{ (\pm 874, 60),
(\pm 1848 , 0),
(\pm 1376 , 0),
(\pm 1036 , 0),
(\pm 899 , 0),
(\pm 849 , 0),
$

$
(\pm 613 , 0),
(\pm 576 , 0),
(\pm 100 , 0),
(-3293 , 0)\}
$
(Figure~\ref{picture_15.png}).

$f = 1748$, $v = 336$, $w = 1780$, $\operatorname{diam(\mathcal{P})} = 5141$,

which gives $d(m, 2m + 15) \leq 5141$.

\begin{figure}[h!]
\center{\includegraphics[width=1\linewidth]{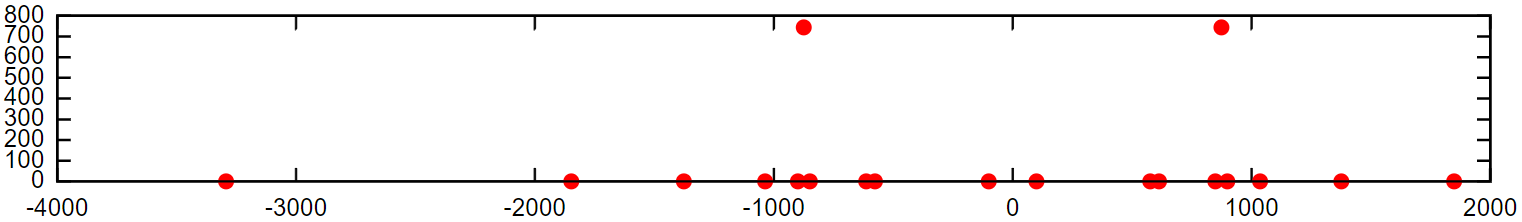}}
\parbox{1\linewidth}{\caption{PIPS of cardinality 19 and diameter 5141}
\label{picture_15.png}}
\end{figure}

\item
$\mathcal{P}=\sqrt{154}/{1} * \{ (\pm 874, 60),
(\pm 3293 , 0),
(\pm 1848 , 0),
(\pm 1376 , 0),
(\pm 1036 , 0),
(\pm 899 , 0),
$

$
(\pm 849 , 0),
(\pm 613 , 0),
(\pm 576 , 0),
(\pm 100 , 0)\}
$
(Figure~\ref{picture_16.png}).

$f = 1748$, $v = 336$, $w = 1780$, $\operatorname{diam(\mathcal{P})} = 6586$,

which gives $d(m, 2m + 16) \leq 6586$.

\begin{figure}[h!]
\center{\includegraphics[width=1\linewidth]{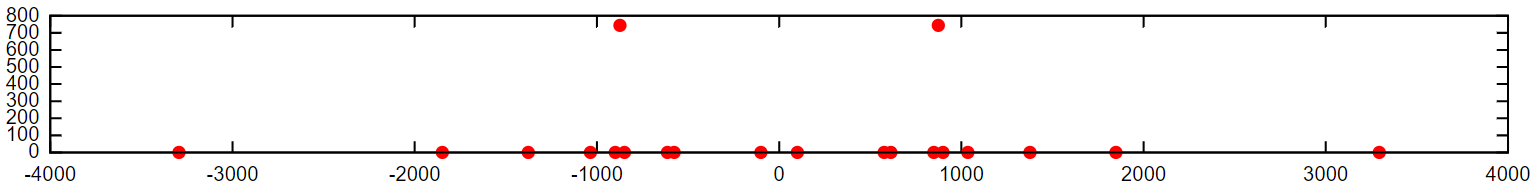}}
\parbox{1\linewidth}{\caption{PIPS of cardinality 20 and diameter 6586}
\label{picture_16.png}}
\end{figure}

\item
$\mathcal{P}=\sqrt{154}/{1} * \{ (\pm 2622, 180),
(\pm 5544 , 0),
(\pm 4128 , 0),
(\pm 3108 , 0),
(\pm 2697 , 0),
$

$
(\pm 2547 , 0),
(\pm 1839 , 0),
(\pm 1728 , 0),
(\pm 904 , 0),
(\pm 300 , 0),
(-6148 , 0)\}
$
(Figure~\ref{picture_17.png}).

$f = 5244$, $v = 1008$, $w = 5340$, $\operatorname{diam(\mathcal{P})} = 11692$,

which gives $d(m, 2m + 17) \leq 11692$.

\item
$\mathcal{P}=\sqrt{154}/{1} * \{ (\pm 2622, 180),
(\pm 6148 , 0),
(\pm 5544 , 0),
(\pm 4128 , 0),
(\pm 3108 , 0),
$

$
(\pm 2697 , 0),
(\pm 2547 , 0),
(\pm 1839 , 0),
(\pm 1728 , 0),
(\pm 904 , 0),
(\pm 300 , 0)\}
$
(Figure~\ref{picture_18.png}).

$f = 5244$, $v = 1008$, $w = 5340$, $\operatorname{diam(\mathcal{P})} = 12296$,

which gives $d(m, 2m + 18) \leq 12296$.

\begin{figure}[h!]
\center{\includegraphics[width=0.75\linewidth]{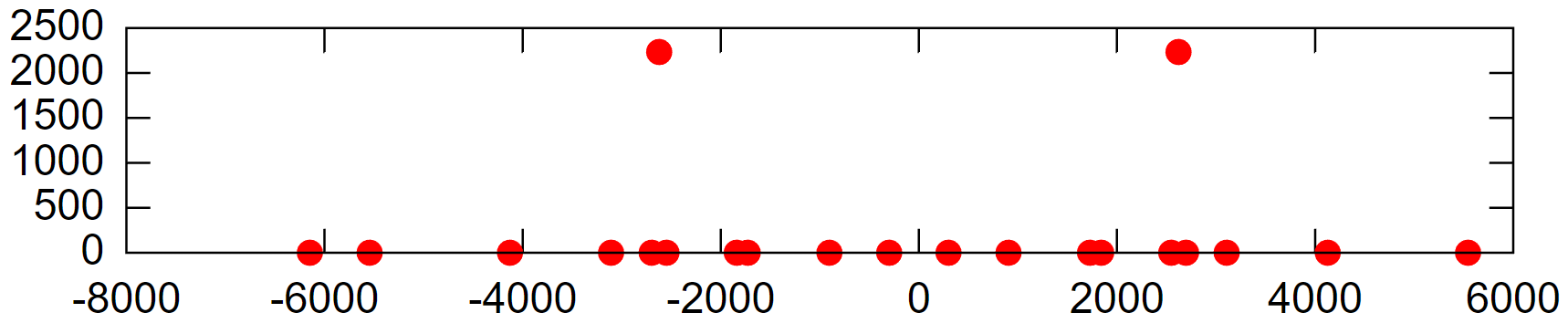}}
\parbox{0.75\linewidth}{\caption{PIPS of cardinality 21 and diameter 11692}
\label{picture_17.png}}
\end{figure}

\item
$\mathcal{P}=\sqrt{154}/{1} * \{ (\pm 2622, 180),
(\pm 6148 , 0),
(\pm 5544 , 0),
(\pm 4128 , 0),
(\pm 3108 , 0),
$

$
(\pm 2697  , 0),
(\pm 2547 , 0),
(\pm 1839 , 0),
(\pm 1728 , 0),
(\pm 904 , 0),
(\pm 300 , 0),
(-9879 , 0)\}
$
(Figure~\ref{picture_19.png}).

$f = 5244$, $v = 1008$, $w = 5340$, $\operatorname{diam(\mathcal{P})} = 16027$,

which gives $d(m, 2m + 19) \leq 16027$.

\begin{figure}[h!]
\center{\includegraphics[width=0.75\linewidth]{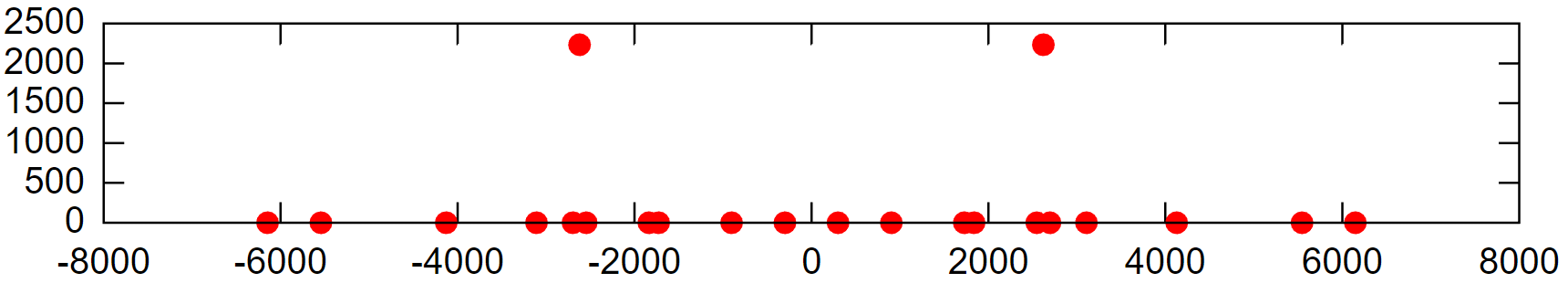}}
\parbox{0.75\linewidth}{\caption{PIPS of cardinality 22 and diameter 12296}
\label{picture_18.png}}
\end{figure}

\item
$\mathcal{P}=\sqrt{154}/{1} * \{ (\pm 2622, 180),
(\pm 9879 , 0),
(\pm 6148 , 0),
(\pm 5544 , 0),
(\pm 4128 , 0),
$

$
(\pm 3108 , 0),
(\pm 2697 , 0),
(\pm 2547 , 0),
(\pm 1839 , 0),
(\pm 1728 , 0),
(\pm 904 , 0),
(\pm 300 , 0)\}
$
(Figure~\ref{picture_20.png}).

$f = 5244$, $v = 1008$, $w = 5340$, $\operatorname{diam(\mathcal{P})} = 19758$,

which gives $d(m, 2m + 20) \leq 19758$.

\begin{figure}[h!]
\center{\includegraphics[width=0.85\linewidth]{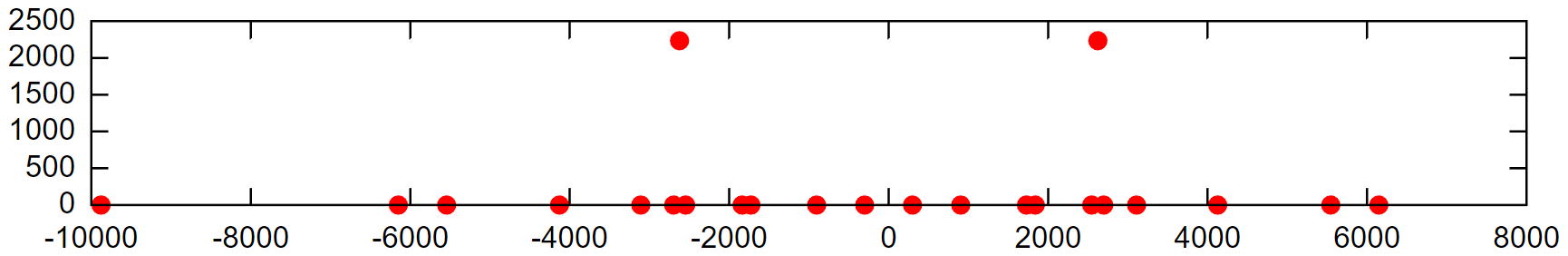}}
\parbox{0.85\linewidth}{\caption{PIPS of cardinality 23 and diameter 16027}
\label{picture_19.png}}
\end{figure}

\item
$\mathcal{P}=\sqrt{154}/{1} * \{ (\pm 5244, 360),
(\pm 12296 , 0),
(\pm 11088 , 0),
(\pm 8579 , 0),
(\pm 8256 , 0),
$

$
(\pm 6216 , 0),
(\pm 5394 , 0),
(\pm 5094 , 0),
(\pm 3678 , 0),
(\pm 3456 , 0),
(\pm 1808 , 0),
(\pm 600 , 0),
$

$
(-19758 , 0)\}
$
(Figure~\ref{picture_21.png}).

$f = 10488$, $v = 1015$, $w = 10537$, $\operatorname{diam(\mathcal{P})} = 32054$,

which gives $d(m, 2m + 21) \leq 32054$.

\begin{figure}[h!]
\center{\includegraphics[width=0.85\linewidth]{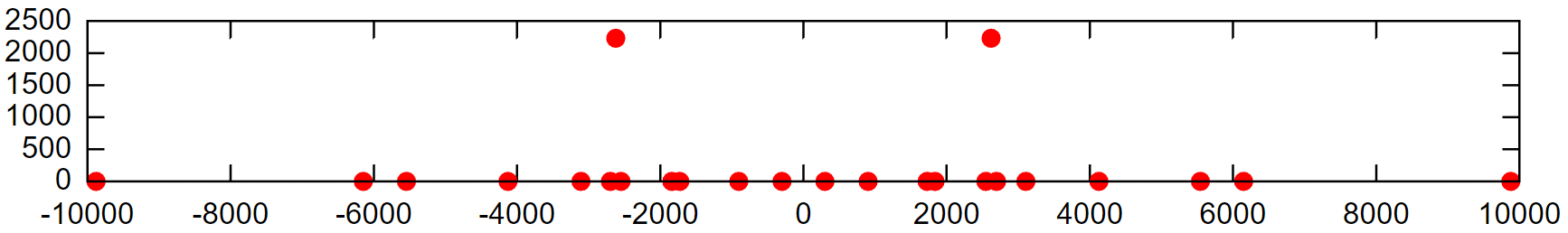}}
\parbox{0.85\linewidth}{\caption{PIPS of cardinality 24 and diameter 19758}
\label{picture_20.png}}
\end{figure}

\item
$\mathcal{P}=\sqrt{154}/{1} * \{ (\pm 5244, 360),
(\pm 19758 , 0),
(\pm 12296 , 0),
(\pm 11088 , 0),
(\pm 8579 , 0),
$

$
(\pm 8256 , 0),
(\pm 6216 , 0),
(\pm 5394 , 0),
(\pm 5094 , 0),
(\pm 3678 , 0),
(\pm 3456 , 0),
(\pm 1808 , 0),
$

$
(\pm 600 , 0)\}
$
(Figure~\ref{picture_22.png}).

$f = 10488$, $v = 1015$, $w = 10537$, $\operatorname{diam(\mathcal{P})} = 39516$,

which gives $d(m, 2m + 22) \leq 39516$.

\begin{figure}[h!]
\center{\includegraphics[width=1\linewidth]{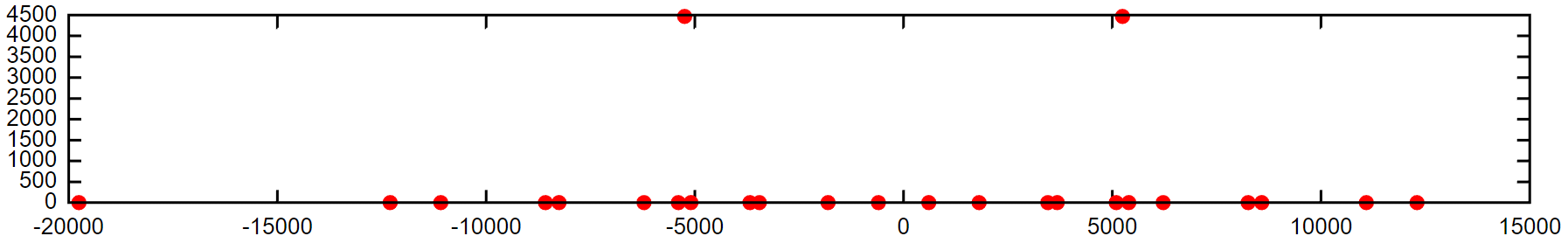}}
\parbox{1\linewidth}{\caption{PIPS of cardinality 25 and diameter 32054}
\label{picture_21.png}}
\end{figure}

\item
$\mathcal{P}=\sqrt{154}/{1} * \{ (\pm 36708, 2520),
(\pm 116058 , 0),
(\pm 86072 , 0),
(\pm 77616 , 0),
(\pm 60053 , 0),
$

$
(\pm 57792 , 0),
(\pm 43512 , 0),
(\pm 37758 , 0),
(\pm 35658 , 0),
(\pm 25746 , 0),
(\pm 24192 , 0),
$

$
(\pm 12656 , 0),
(\pm 4200 , 0),
(-138306 , 0)\}
$
(Figure~\ref{picture_23.png}).

$f = 73416$, $v = 2710$, $w = 73466$, $\operatorname{diam(\mathcal{P})} = 254364$,

which gives $d(m, 2m + 23) \leq 254364$.

\begin{figure}[h!]
\center{\includegraphics[width=1\linewidth]{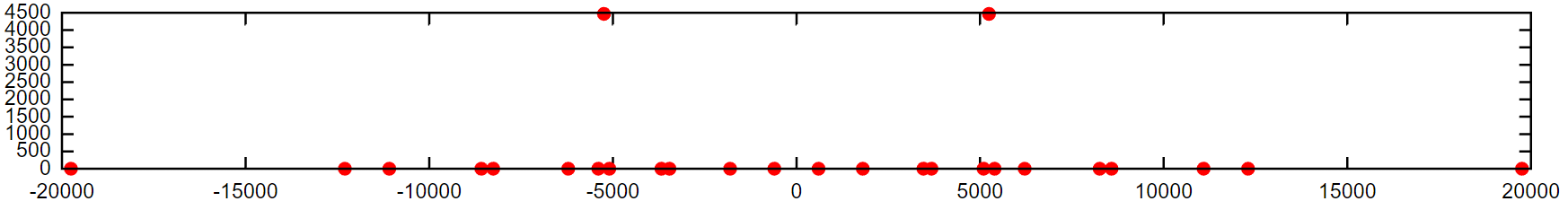}}
\parbox{1\linewidth}{\caption{PIPS of cardinality 26 and diameter 39516}
\label{picture_22.png}}
\end{figure}

\item
$\mathcal{P}=\sqrt{154}/{1} * \{ (\pm 36708, 2520),
(\pm 138306 , 0),
(\pm 116058 , 0),
(\pm 86072 , 0),
(\pm 77616 , 0),
$

$
(\pm 60053 , 0),
(\pm 57792 , 0),
(\pm 43512 , 0),
(\pm 37758 , 0),
(\pm 35658 , 0),
(\pm 25746 , 0),
$

$
(\pm 24192 , 0),
(\pm 12656 , 0),
(\pm 4200 , 0)\}
$
(Figure~\ref{picture_24.png}).

$f = 73416$, $v = 2710$, $w = 73466$, $\operatorname{diam(\mathcal{P})} = 276612$,

which gives $d(m, 2m + 24) \leq 276612$.

\begin{figure}[h!]
\center{\includegraphics[width=1\linewidth]{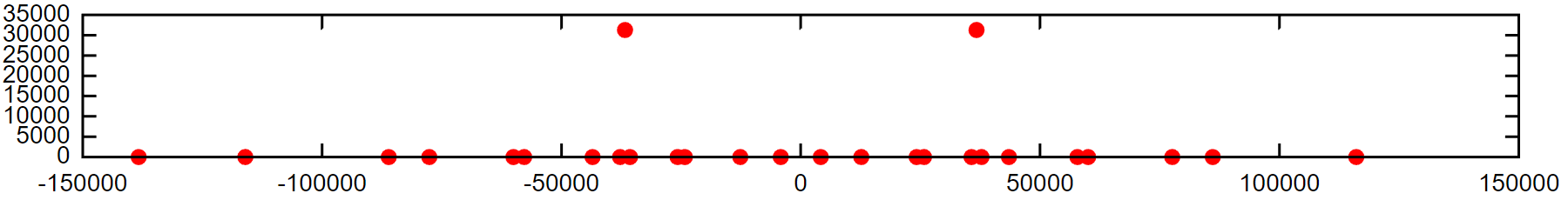}}
\parbox{1\linewidth}{\caption{PIPS of cardinality 27 and diameter 254364}
\label{picture_23.png}}
\end{figure}

\item
$\mathcal{P}=\sqrt{154}/{1} * \{ (\pm 36708, 2520),
(\pm 138306 , 0),
(\pm 116058 , 0),
(\pm 86072 , 0),
(\pm 77616 , 0),
$

$
(\pm 60053 , 0),
(\pm 57792 , 0),
(\pm 43512 , 0),
(\pm 37758 , 0),
(\pm 35658 , 0),
(\pm 25746 , 0),
$

$
(\pm 24192 , 0),
(\pm 12656 , 0),
(\pm 4200 , 0),
(-199863 , 0)\}
$
(Figure~\ref{picture_25.png}).

$f = 73416$, $v = 2710$, $w = 73466$, $\operatorname{diam(\mathcal{P})} = 338169$,

which gives $d(m, 2m + 25) \leq 338169$.

\begin{figure}[h!]
\center{\includegraphics[width=1\linewidth]{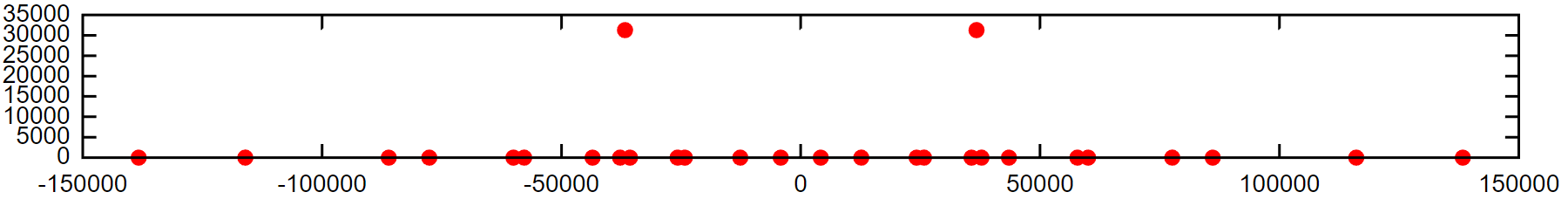}}
\parbox{1\linewidth}{\caption{PIPS of cardinality 28 and diameter 276612}
\label{picture_24.png}}
\end{figure}

\begin{figure}[h!]
\center{\includegraphics[width=1\linewidth]{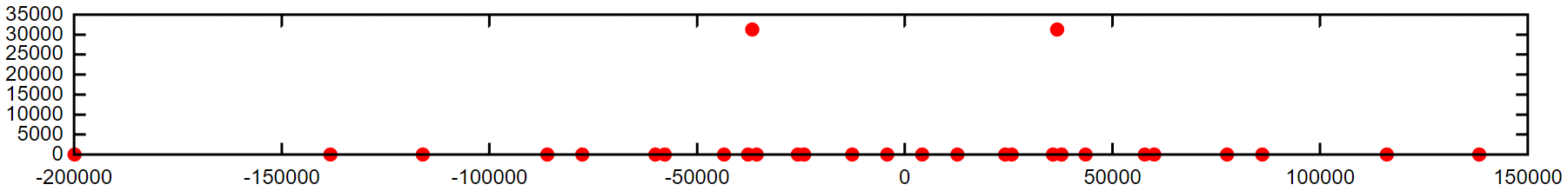}}
\parbox{1\linewidth}{\caption{PIPS of cardinality 29 and diameter 338169}
\label{picture_25.png}}
\end{figure}

\item
$\mathcal{P}=\sqrt{154}/{1} * \{ (\pm 36708, 2520),
(\pm 199863 , 0),
(\pm 138306 , 0),
(\pm 116058 , 0),
(\pm 86072 , 0),
$

$
(\pm 77616 , 0),
(\pm 60053 , 0),
(\pm 57792 , 0),
(\pm 43512 , 0),
(\pm 37758 , 0),
(\pm 35658 , 0),
$

$
(\pm 25746 , 0),
(\pm 24192 , 0),
(\pm 12656 , 0),
(\pm 4200 , 0)\}
$
(Figure~\ref{picture_26.png}).

$f = 73416$, $v = 2710$, $w = 73466$, $\operatorname{diam(\mathcal{P})} = 399726$,

which gives $d(m, 2m + 26) \leq 399726$.

\begin{figure}[h!]
\center{\includegraphics[width=1\linewidth]{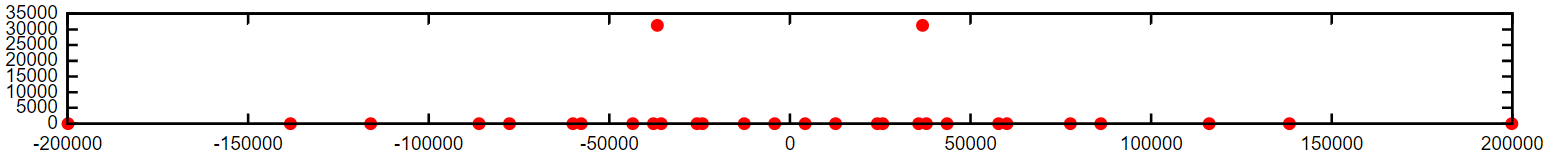}}
\parbox{1\linewidth}{\caption{PIPS of cardinality 30 and diameter 399726}
\label{picture_26.png}}
\end{figure}

\item
$\mathcal{P}=\sqrt{154}/{1} * \{ (\pm 580290, 0),
(\pm 430360, 0),
(\pm 388080, 0),
(\pm 300265, 0),
$

$
(\pm 288960, 0),
(\pm 217560, 0),
(\pm 188790, 0),
(\pm 183540, 12600),
(\pm 178290, 0),
$

$
(\pm 128730, 0),
(\pm 120960, 0),
(\pm 63280, ),
(\pm 21000, 0),
(\pm 11224, 0),
$

$
(-1033912, 0),
(-999315, 0),
(-691530  , 0)\}
$
(Figure~\ref{picture_31.png}).

$f = 367080$, $v = 3206$, $w = 367094$, $\operatorname{diam(\mathcal{P})} = 1614202$,

which gives $d(m, 2m + 27) \leq 1614202$.

\begin{figure}[h!]
\center{\includegraphics[width=1\linewidth]{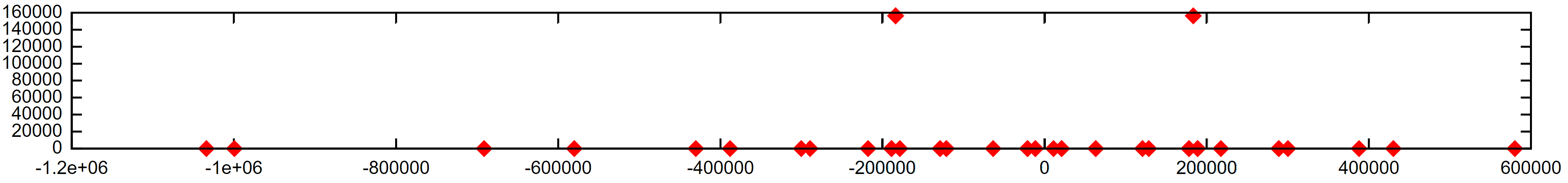}}
\parbox{1\linewidth}{\caption{PIPS of cardinality 31 and diameter 1614202}
\label{picture_31.png}}
\end{figure}

\item
$\mathcal{P}=\sqrt{154}/{1} * \{ (\pm 691530, 0),
(\pm 580290, 0),
(\pm 430360, 0),
(\pm 388080, 0),
$

$
(\pm 300265, 0),
(\pm 288960, 0),
(\pm 217560, 0),
(\pm 188790, 0),
(\pm 183540, 12600),
$

$
(\pm 178290, 0),
(\pm 128730, 0),
(\pm 120960, 0),
(\pm 63280, 0),
(\pm 21000, 0),
$

$
(\pm 11224, 0),
(-1033912, 0),
(-999315, 0)\}
$
(Figure~\ref{picture_32.png}).

$f = 367080$, $v = 3206$, $w = 367094$, $\operatorname{diam(\mathcal{P})} = 1725442$,

which gives $d(m, 2m + 28) \leq 1725442$.

\begin{figure}[h!]
\center{\includegraphics[width=1\linewidth]{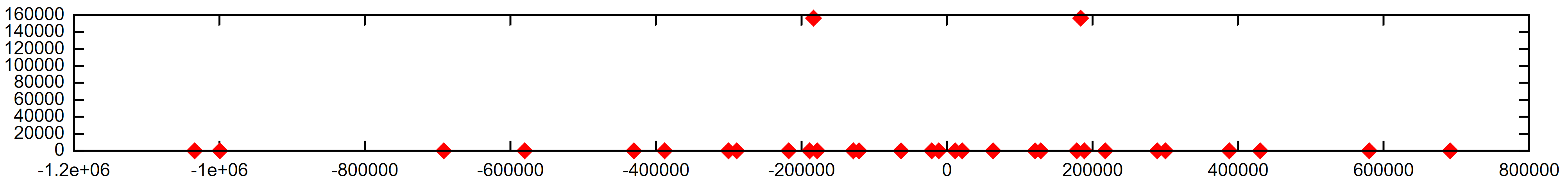}}
\parbox{1\linewidth}{\caption{PIPS of cardinality 32 and diameter 1725442}
\label{picture_32.png}}
\end{figure}

\item
$\mathcal{P}=\sqrt{154}/{1} * \{ (\pm 999315, 0),
(\pm 691530, 0),
(\pm 580290, 0),
(\pm 430360, 0),
$

$
(\pm 388080, 0),
(\pm 300265, 0),
(\pm 288960, 0),
(\pm 217560, 0),
(\pm 188790, 0),
$

$
(\pm 183540, 12600),
(\pm 178290, 0),
(\pm 128730, 0),
(\pm 120960, 0),
(\pm 63280, 0),
$

$
(\pm 21000, 0),
(\pm 11224, 0),
(-1033912, 0)\}
$
(Figure~\ref{picture_33.png}).

$f = 367080$, $v = 3206$, $w = 367094$, $\operatorname{diam(\mathcal{P})} = 2033227$,

which gives $d(m, 2m + 29) \leq 2033227$.

\begin{figure}[h!]
\center{\includegraphics[width=1\linewidth]{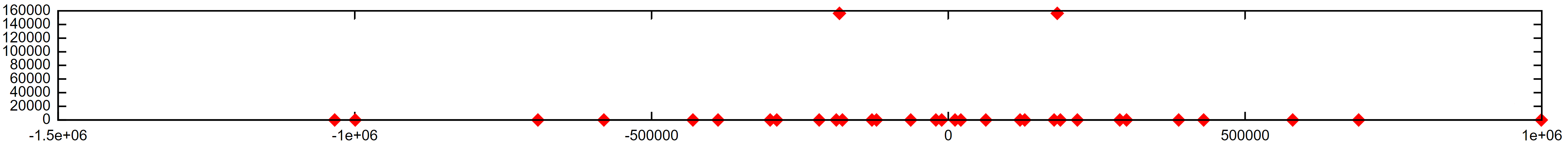}}
\parbox{1\linewidth}{\caption{PIPS of cardinality 33 and diameter 2033227}
\label{picture_33.png}}
\end{figure}

\item
$\mathcal{P}=\sqrt{154}/{1} * \{ (\pm 1033912, 0),
(\pm 999315, 0),
(\pm 691530, 0),
(\pm 580290, 0),
$

$
(\pm 430360, 0),
(\pm 388080, 0),
(\pm 300265, 0),
(\pm 288960, 0),
(\pm 217560, 0),
$

$
(\pm 188790, 0),
(\pm 183540, 12600),
(\pm 178290, 0),
(\pm 128730, 0),
(\pm 120960, 0),
$

$
(\pm 63280, 0),
(\pm 21000, 0),
(\pm 11224, 0)\}
$
(Figure~\ref{picture_34.png}).

$f = 367080$, $v = 3206$, $w = 367094$, $\operatorname{diam(\mathcal{P})} = 2067824$,

which gives $d(m, 2m + 30) \leq 2067824$.

\begin{figure}[h!]
\center{\includegraphics[width=1\linewidth]{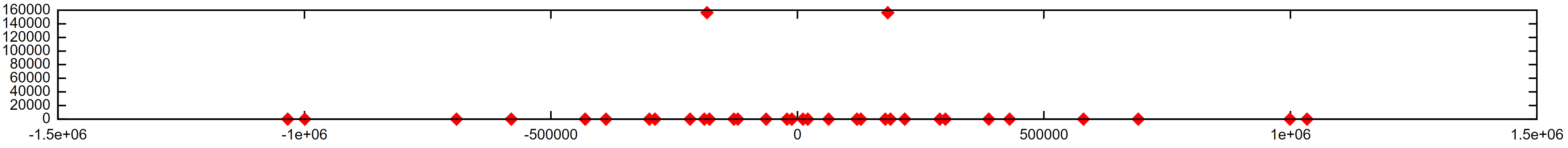}}
\parbox{1\linewidth}{\caption{PIPS of cardinality 34 and diameter 2067824}
\label{picture_34.png}}
\end{figure}

\item
$\mathcal{P}=\sqrt{154}/{1} * \{ (\pm 3997260, 0),
(\pm 2766120, 0),
(\pm 2321160, 0),
(\pm 1721440, 0),
$

$
(\pm 1552320, 0),
(\pm 1201060, 0),
(\pm 1155840, 0),
(\pm 870240, 0),
(\pm 755160, 0),
$

$
(\pm 734160, 50400),
(\pm 713160, 0),
(\pm 514920, 0),
(\pm 483840, 0),
(\pm 308175, 0),
$

$
(\pm 253120, 0),
(\pm 84000, 0),
(\pm 44896, 0),
(-4135648 , 0)\}
$
(Figure~\ref{picture_35.png}).

$f = 1468320$, $v = 12824$, $w = 1468376$, $\operatorname{diam(\mathcal{P})} = 8132908$,

which gives $d(m, 2m + 31) \leq 8132908$.

\begin{figure}[h!]
\center{\includegraphics[width=1\linewidth]{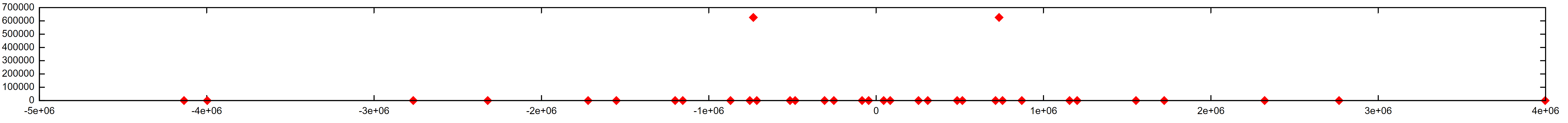}}
\parbox{1\linewidth}{\caption{PIPS of cardinality 35 and diameter 8132908}
\label{picture_35.png}}
\end{figure}

\item
$\mathcal{P}=\sqrt{154}/{1} * \{ (\pm 4135648, 0),
(\pm 3997260, 0),
(\pm 2766120, 0),
(\pm 2321160, 0),
$

$
(\pm 1721440, 0),
(\pm 1552320, 0),
(\pm 1201060, 0),
(\pm 1155840, 0),
(\pm 870240, 0),
$

$
(\pm 755160, 0),
(\pm 734160, 50400),
(\pm 713160, 0),
(\pm 514920, 0),
(\pm 483840, 0),
$

$
(\pm 308175, 0),
(\pm 253120, 0),
(\pm 84000, 0),
(\pm 44896, 0)\}
$
(Figure~\ref{picture_36.png}).

$f = 1468320$, $v = 12824$, $w = 1468376$, $\operatorname{diam(\mathcal{P})} = 8271296$,

which gives $d(m, 2m + 32) \leq 8271296$.

\begin{figure}[h!]
\center{\includegraphics[width=1\linewidth]{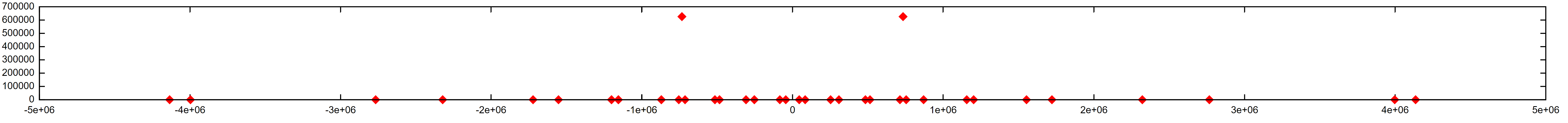}}
\parbox{1\linewidth}{\caption{PIPS of cardinality 36 and diameter 8271296}
\label{picture_36.png}}
\end{figure}

\item
$\mathcal{P}=\sqrt{154}/{1} * \{ (\pm 12991095, 0),
(\pm 8989890, 0),
(\pm 7543770, 0),
(\pm 5594680, 0),
$

$
(\pm 5045040, 0),
(\pm 3903445, 0),
(\pm 3756480, 0),
(\pm 3694845, 0),
(\pm 2828280, 0),
$

$
(\pm 2454270, 0),
(\pm 2386020, 163800),
(\pm 2317770, 0),
(\pm 1673490, 0),
(\pm 1572480, 0),
$

$
(\pm 822640, 0),
(\pm 484680, 0),
(\pm 273000, 0),
(\pm 145912, 0),
$

$
(-13440856, 0)\}
$
(Figure~\ref{picture_37.png}).

$f = 4772040$, $v = 41678$, $w = 4772222$, $\operatorname{diam(\mathcal{P})} = 26431951$,

which gives $d(m, 2m + 33) \leq 26431951$.

\begin{figure}[h!]
\center{\includegraphics[width=1\linewidth]{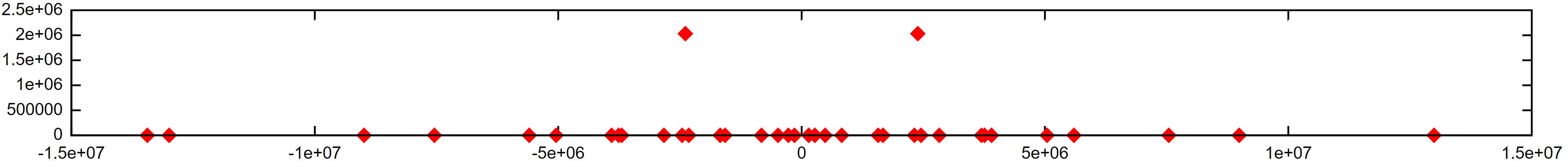}}
\parbox{1\linewidth}{\caption{PIPS of cardinality 37 and diameter 26431951}
\label{picture_37.png}}
\end{figure}

\item
$\mathcal{P}=\sqrt{154}/{1} * \{ (\pm 13440856, 0),
(\pm 12991095, 0),
(\pm 8989890, 0),
(\pm 7543770, 0),
$

$
(\pm 5594680, 0),
(\pm 5045040, 0),
(\pm 3903445, 0),
(\pm 3756480, 0),
(\pm 3694845, 0),
$

$
(\pm 2828280, 0),
(\pm 2454270, 0),
(\pm 2386020, 163800),
(\pm 2317770, 0),
(\pm 1673490, 0),
$

$
(\pm 1572480, 0),
(\pm 822640, 0),
(\pm 484680, 0),
(\pm 273000, 0),
$

$
(\pm 145912, 0)\}
$
(Figure~\ref{picture_38.png}).

$f = 4772040$, $v = 41678$, $w = 4772222$, $\operatorname{diam(\mathcal{P})} = 26881712$,

which gives $d(m, 2m + 34) \leq 26881712$.

\begin{figure}[h!]
\center{\includegraphics[width=1\linewidth]{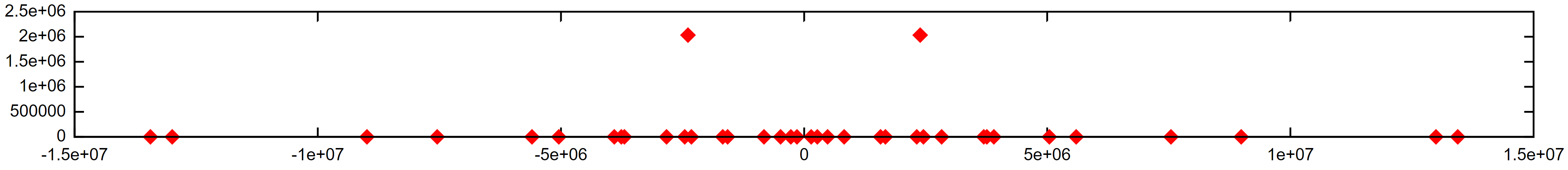}}
\parbox{1\linewidth}{\caption{PIPS of cardinality 38 and diameter 26881712}
\label{picture_38.png}}
\end{figure}

\item
$\mathcal{P}=\sqrt{154}/{1} * \{ (\pm 51964380, 0),
(\pm 35959560, 0),
(\pm 30175080, 0),
(\pm 22378720, 0),
$

$
(\pm 20180160, 0),
(\pm 15613780, 0),
(\pm 15025920, 0),
(\pm 14779380, 0),
(\pm 11313120, 0),
$

$
(\pm 9817080, 0),
(\pm 9544080, 655200),
(\pm 9271080, 0),
(\pm 6693960, 0),
(\pm 6289920, 0),
$

$
(\pm 4006275, 0),
(\pm 3290560, 0),
(\pm 1938720, 0),
(\pm 1092000, 0),
$

$
(\pm 583648, 0),
(-53763424 , 0)\}
$
(Figure~\ref{picture_39.png}).

$f = 19088160$, $v = 8738$, $w = 19088162$, $\operatorname{diam(\mathcal{P})} = 105727804$,

which gives $d(m, 2m + 35) \leq 105727804$.

\begin{figure}[h!]
\center{\includegraphics[width=1\linewidth]{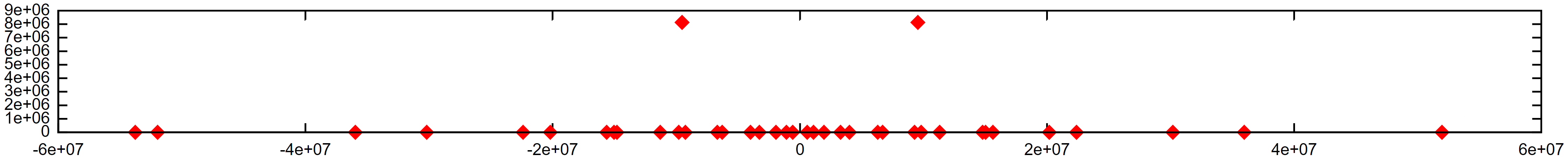}}
\parbox{1\linewidth}{\caption{PIPS of cardinality 39 and diameter 105727804}
\label{picture_39.png}}
\end{figure}

\begin{figure}[h!]
	\begin{center}
	\includegraphics[width=1\linewidth]{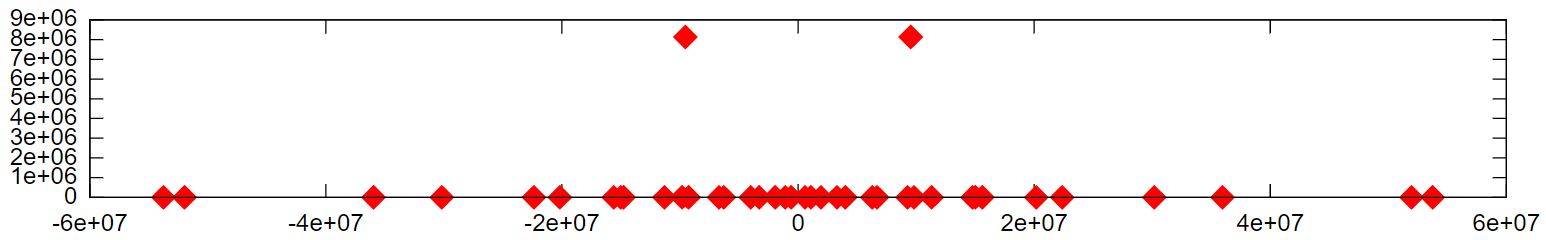}
	\parbox{1\linewidth}{\caption{PIPS of cardinality 40 and diameter 107526848}
	\label{picture_40.png}}
	\end{center}
\end{figure}

\item
$\mathcal{P}=\sqrt{154}/{1} * \{ (\pm 53763424, 0),
(\pm 51964380 , 0),
(\pm 35959560 , 0),
(\pm 30175080 , 0),
$

$
(\pm 22378720 , 0),
(\pm 20180160 , 0),
(\pm 15613780 , 0),
(\pm 15025920 , 0),
(\pm 14779380 , 0),
$

$
(\pm 11313120 , 0),
(\pm 9817080 , 0),
(\pm 9544080 , 655200),
(\pm 9271080 , 0),
(\pm 6693960 , 0),
$

$
(\pm 6289920 , 0),
(\pm 4006275 , 0),
(\pm 3290560 , 0),
(\pm 1938720 , 0),
(\pm 1092000 , 0),
$

$
(\pm 583648 , 0)\}
$
(Figure~\ref{picture_40.png}).

$f = 19088160$, $v = 8738$, $w = 19088162$, $\operatorname{diam}(\mathcal{P})
= 107526848$,

which gives $d(m, 2m + 36) \leq 107526848$.

\end{itemize}

\begin{remark}
In all the examples above, the maximum in the right-hand side of (\ref{formula1})
is attained on the diameter of the system.
The following example shows that the maximum can be attained on $w$:
\end{remark}

\begin{figure}[h!]
	\includegraphics[width=.48\linewidth]{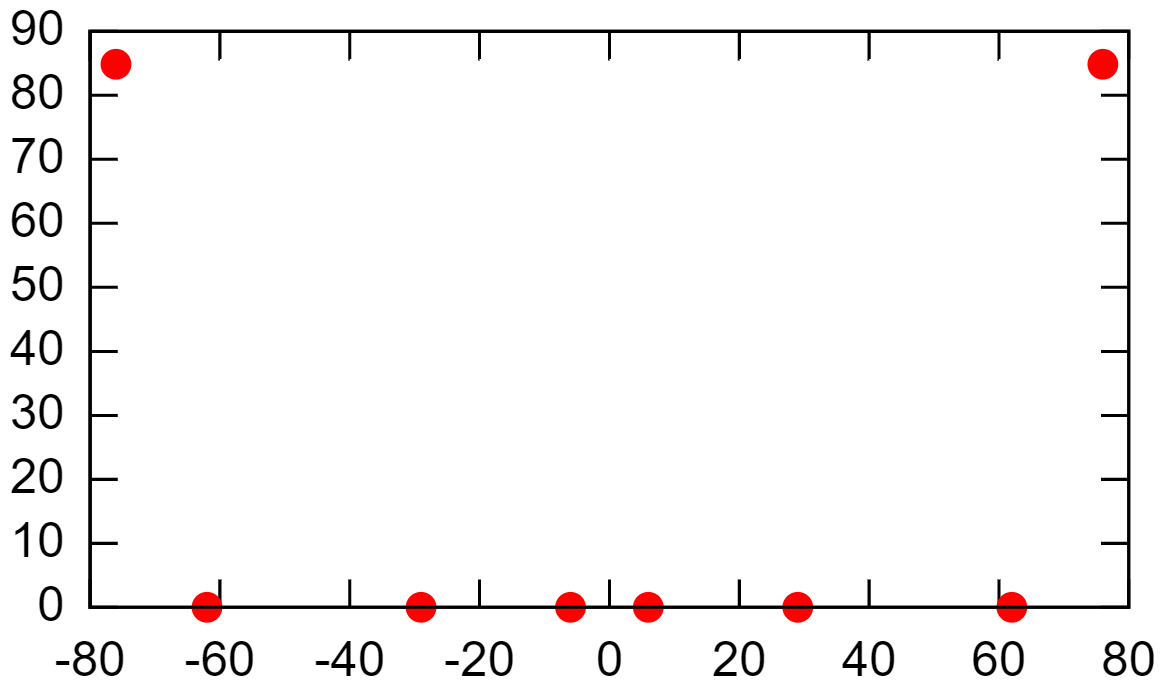}
	\hfill
	\includegraphics[width=.48\linewidth]{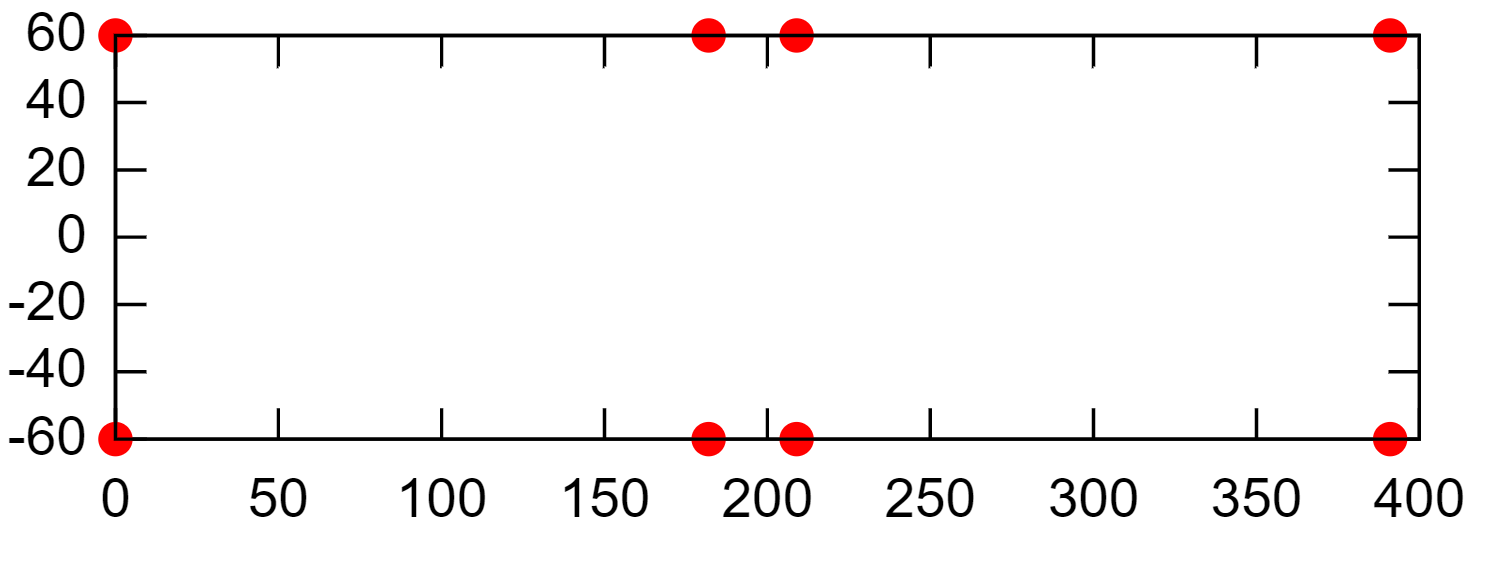}
	\\
	\parbox{.48\linewidth}{\caption{PIPS of cardinality 8 and diameter 162}
	\label{picture_10.png}}
	\hfill
	\parbox{.48\linewidth}{\caption{PIPS of cardinality 8 and diameter 409}
	\label{picture.png}}
\end{figure}

\begin{itemize}
\setlength{\itemsep}{-1mm}

\item
$\mathcal{P}=\sqrt{2}/{1} * \{ (\pm 76, 60),
(\pm 62 , 0),
(\pm 29 , 0),
(\pm 6 , 0)\}
$
(Figure~\ref{picture_10.png}).

$f = 152$, $v = 114$, $w = 190$, $\operatorname{diam(\mathcal{P})} = 162$.

We can obtain the estimate $d(m, 2m + 2) \leq 190$ (which
is weaker than the one given in~\eqref{eq:d_m_2m+2}.

\end{itemize}

Guy gives \cite[D 20]{guy2013unsolved} a PIPS of $8$ points located
on two parallel lines. The points of the set form a rectangle
(Figure~\ref{rectangle.pdf}):

\begin{figure}[htbp]
	\includegraphics[width=.6\linewidth]{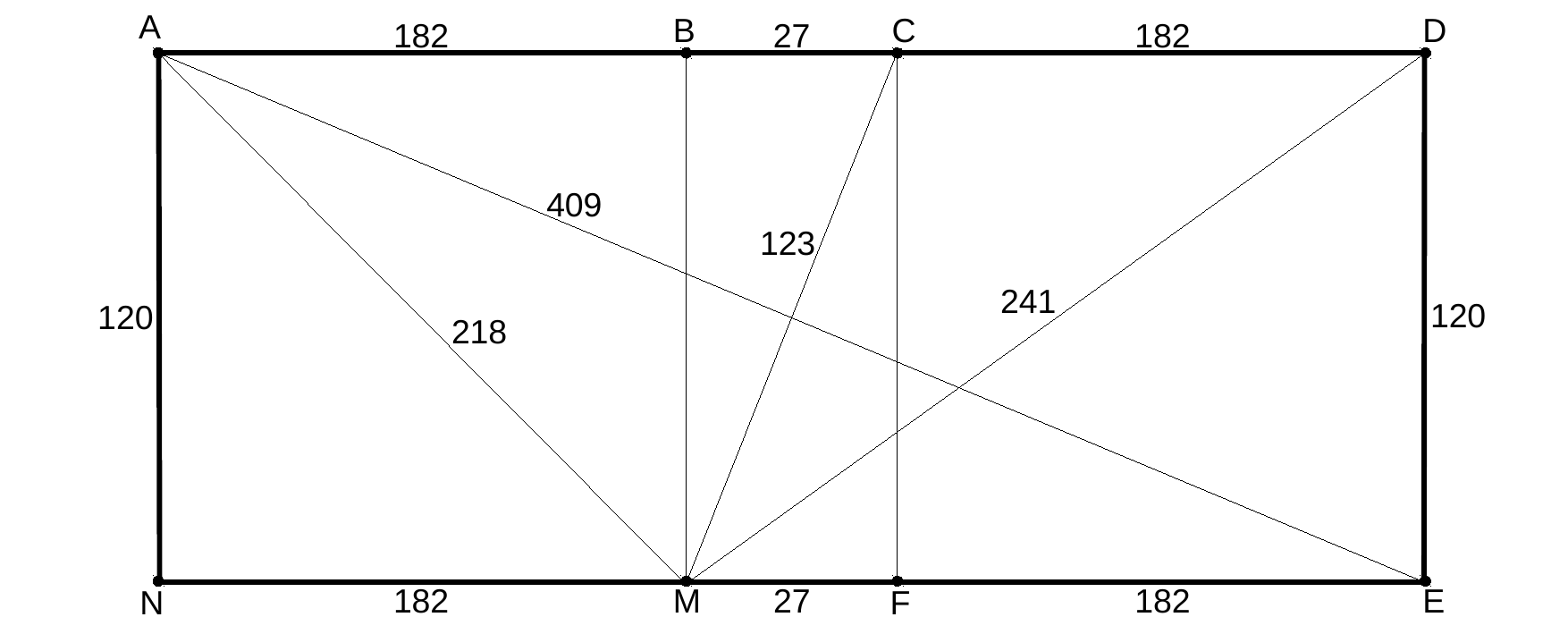}
	\hfill
	\includegraphics[width=.5\linewidth]{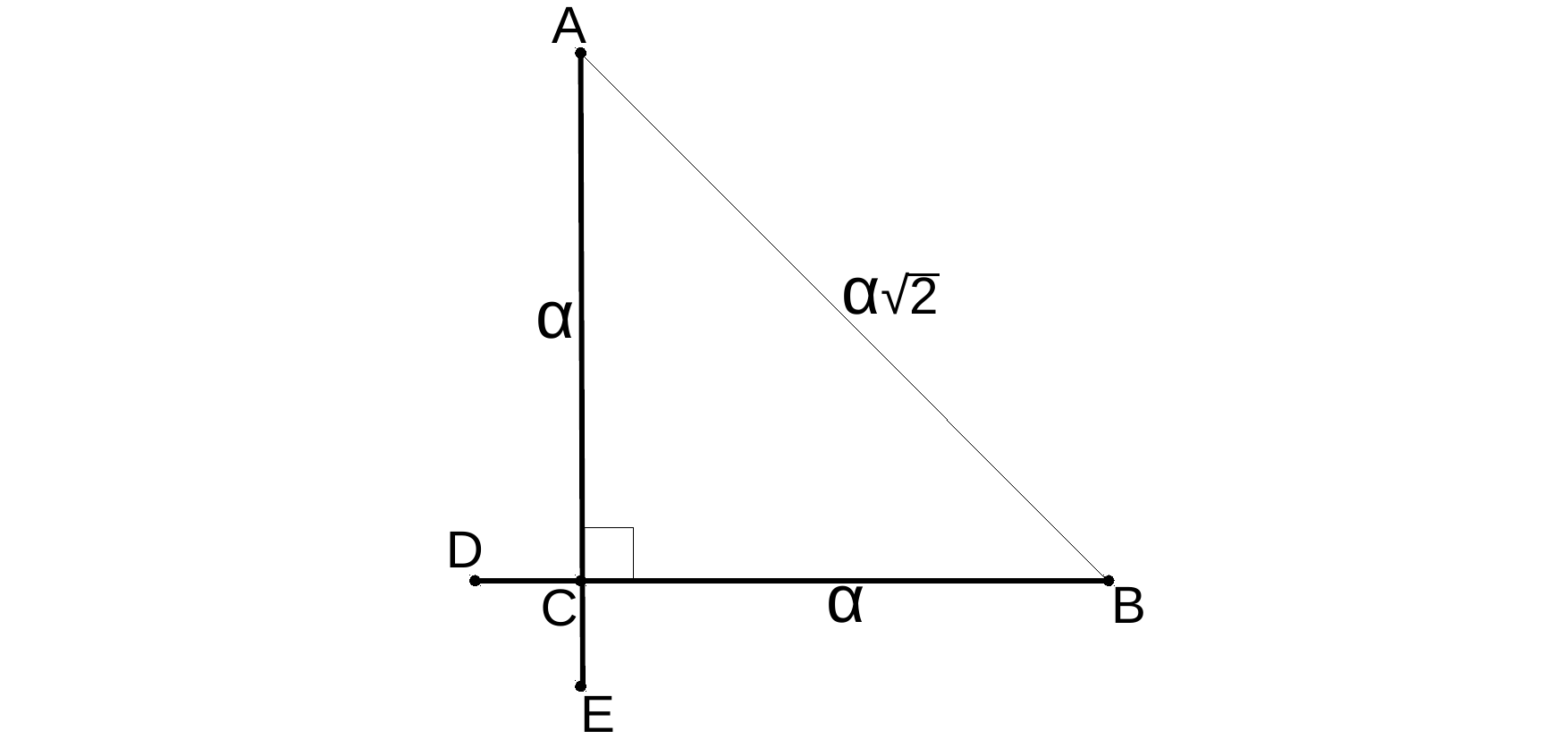}
	\\
	\parbox{.65\linewidth}{\caption{Distances}
	\label{rectangle.pdf}}
	\hfill
	\parbox{.5\linewidth}{\caption{Contradiction}
	\label{right_angle.pdf}}
\end{figure}

\begin{itemize}
\setlength{\itemsep}{-1mm}

\item
$\mathcal{P}=\sqrt{1}/{1} * \{ (0, \pm 60),
(182 , \pm 60),
(209 , \pm 60),
(391 , \pm 60)\}
$
(Figure~\ref{picture.png}),

where $n = 8$, $\operatorname{diam(\mathcal{P})} = 409$. Using the ``blowing up''
construction, we obtain
\begin{equation}\label{result1}
d(m, 4m) \leq 409
\end{equation}
Figure~\ref{picture_11.pdf} shows an example for $m = 3$.

\begin{figure}[h!]
\center{\includegraphics[width=1\linewidth]{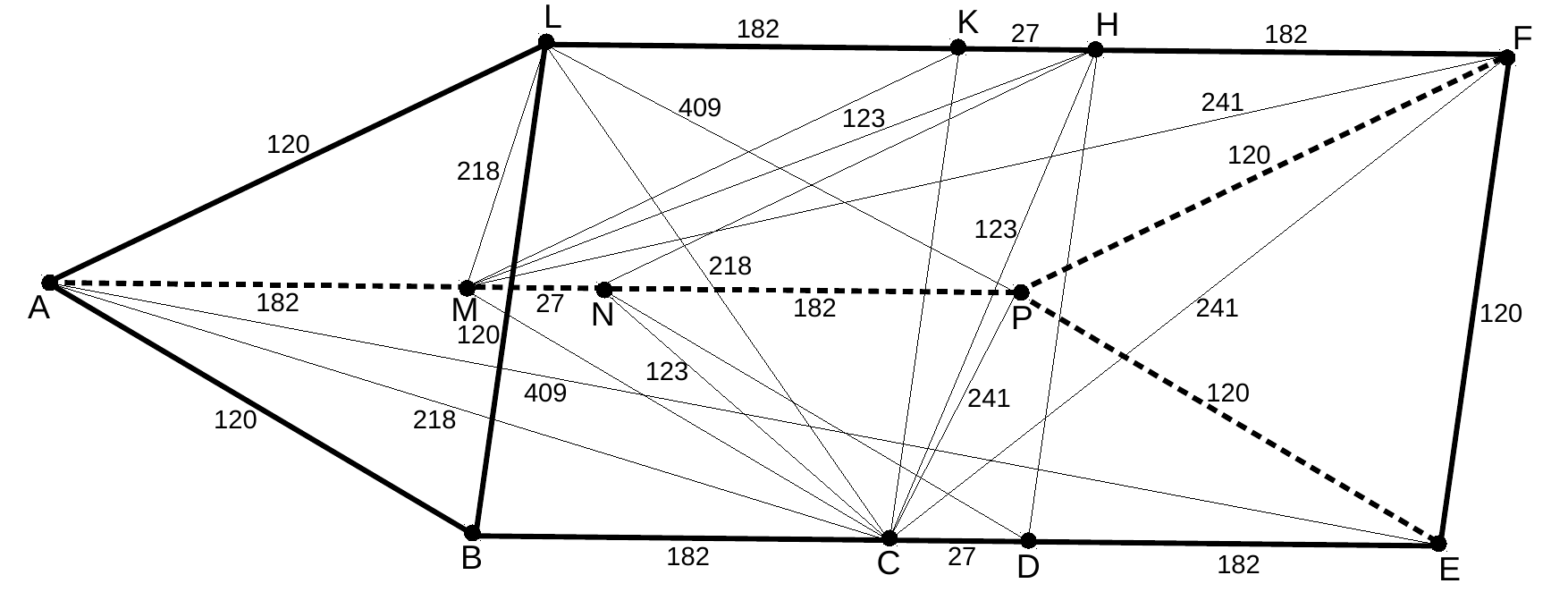}}
\parbox{1\linewidth}{\caption{IPS of cardinality 12 and diameter 409}
\label{picture_11.pdf}}
\end{figure}

\end{itemize}

\section{Bounds based on pyramid sets}
\label{sec:bounds_pyramid}

\begin{theorem}
Let $\mathcal{P}$ be a planar integral point set consisting of $k$
points on the line $l_{1}$ and $k$ points on the line $l_{2}$. Besides, these
points are symmetric with respect to one of the bisectors of the angles
formed by the intersection of lines $l_{1}$ and $l_{2}$, then

\begin{equation}\label{formula2}
d(m, (k - \alpha)m + \alpha) \leq \operatorname{diam}(\mathcal{P}),
\end{equation}

where
\begin{equation*}
\alpha =
\begin{cases}
1, \text{when the intersection point} \in \mathcal{P}, \\
0, \text{when the intersection point} \notin \mathcal{P}.
\end{cases}
\end{equation*}

\end{theorem}

\begin{remark}
The angles can be acute or obtuse. Figure~\ref{right_angle.pdf} shows that the intersection angle of the lines $l_{1}$ and $l_{2}$ cannot be equal to ${\pi}/2$.
\end{remark}

Below we give the examples of planar integral point sets and the corresponding
estimates of the function $d(m, n)$ for $n = 3m + 1$ and $n = 4m + 1$.

\begin{figure}[htbp]
	\includegraphics[width=.48\linewidth]{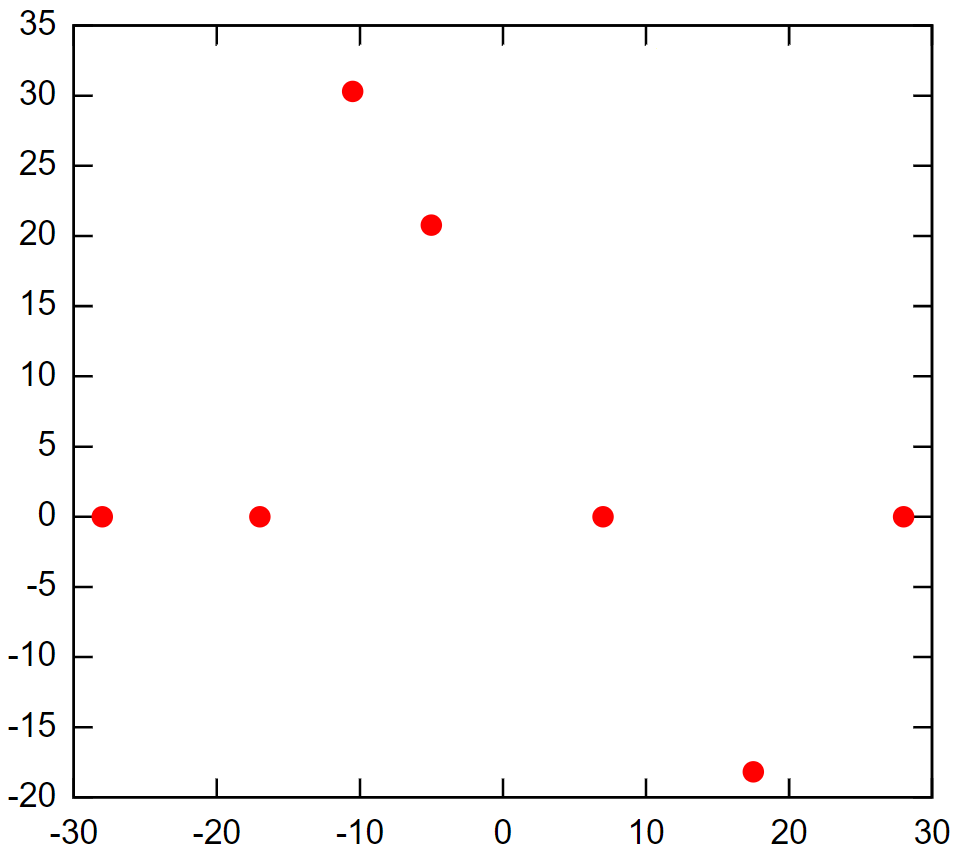}
	\hfill
	\includegraphics[width=.48\linewidth]{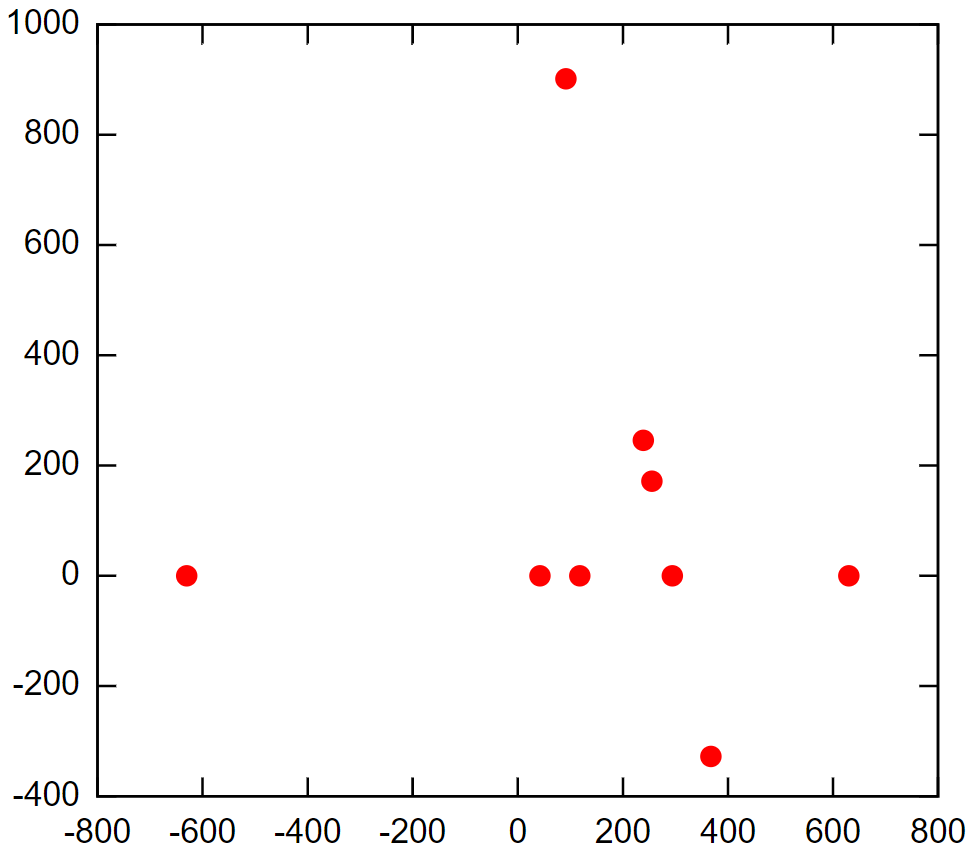}
	\\
	\parbox{.48\linewidth}{\caption{PIPS of cardinality 7 and diameter 56}
	\label{picture_56.png}}
	\hfill
	\parbox{.48\linewidth}{\caption{PIPS of cardinality 9 and diameter 1260}
	\label{picture_1260.png}}
\end{figure}

\begin{figure}[h!]
\center{\includegraphics[width=1\linewidth]{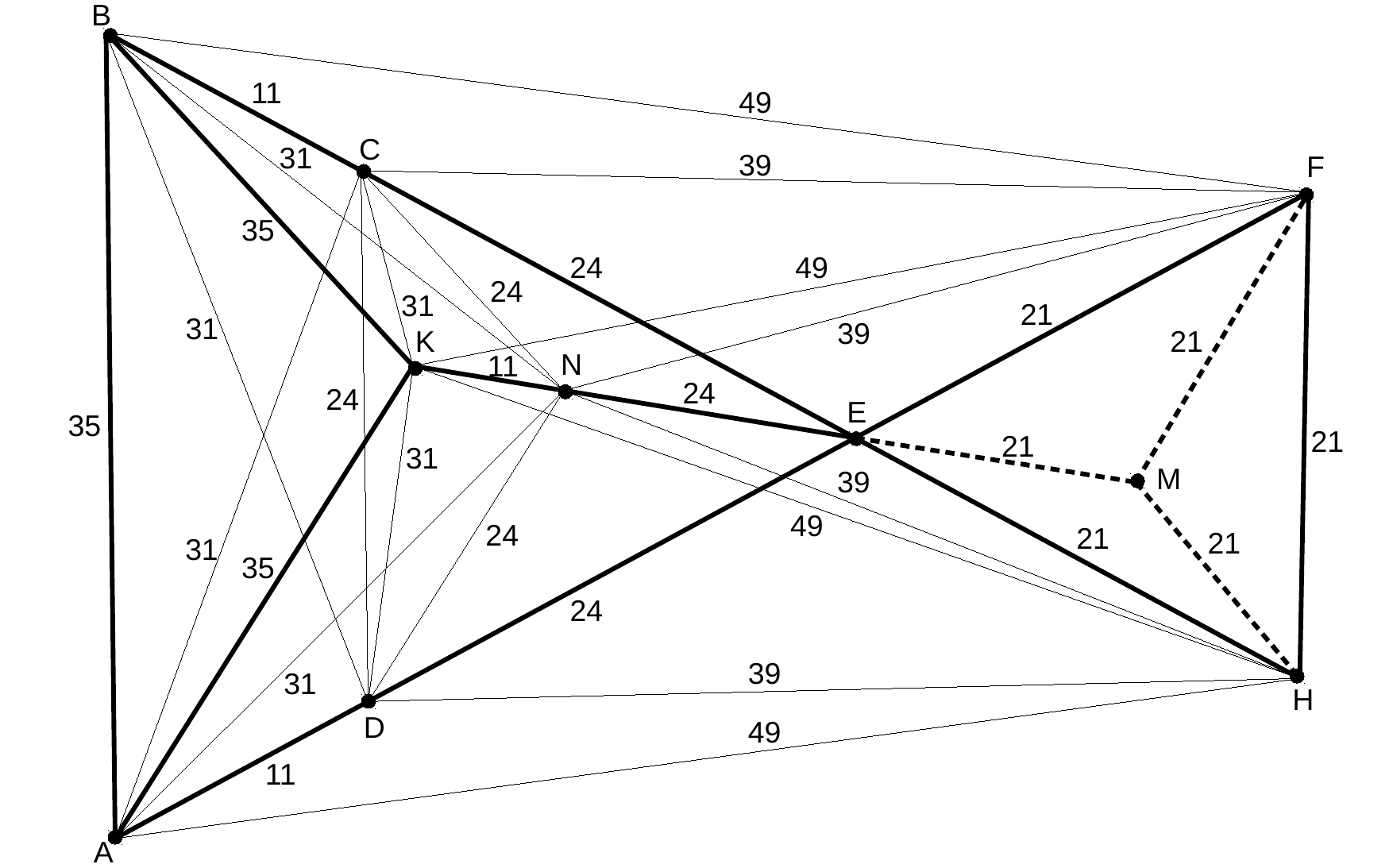}}
\parbox{1\linewidth}{\caption{IPS of cardinality 10 and diameter 56}
\label{picture_12.pdf}}
\end{figure}

\begin{figure}[h!]
\center{\includegraphics[width=1\linewidth]{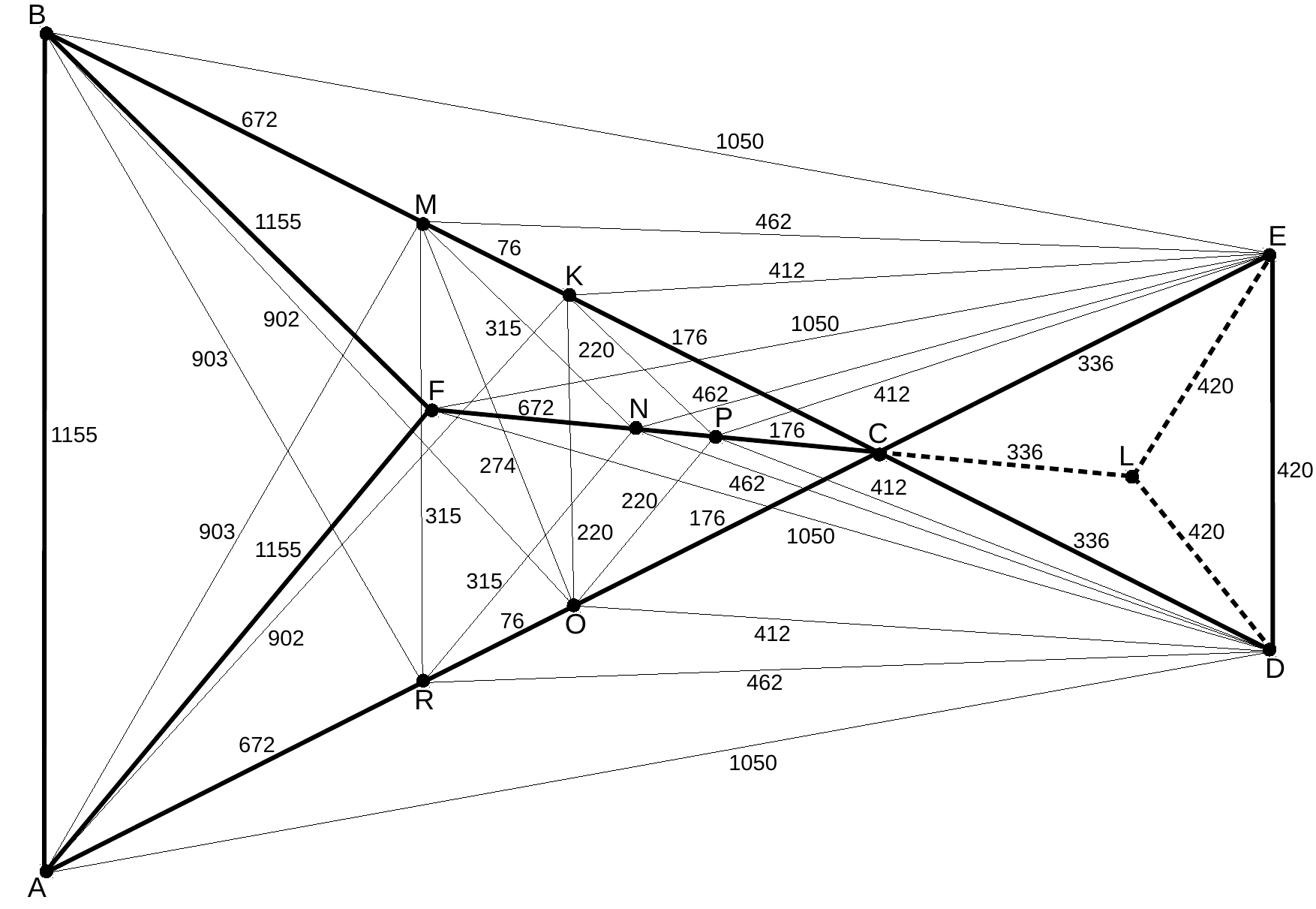}}
\parbox{1\linewidth}{\caption{IPS of cardinality 13 and diameter 1260}
\label{picture_1260_R3.pdf}}
\end{figure}

\begin{itemize}
\setlength{\itemsep}{-1mm}

\item
$\mathcal{P}=\sqrt{3}/{2} * \{ (\pm 56, 0),
(14, 0),
(-34, 0),
(-10, 24),
(-21 , 35),
(35, -21)\}
$
(Figure~\ref{picture_56.png}),

where $n = 7$, $\operatorname{diam(\mathcal{P})} = 56$. Using the ``blowing up''
construction, we obtain
\begin{equation}\label{result2}
d(m, 3m + 1) \leq 56
\end{equation}
The resulting estimate improves the estimate for $d(m, 3m)$, which is presented
in \cite{kemnitz1988punktmengen}. Figure~\ref{picture_12.pdf} shows an example
for $m = 3$.

\item
$\mathcal{P}=\sqrt{39}/{8} * \{ (\pm 5040, 0),
(1911, 315),
(2352, 0),
(944, 0),
(336, 0),
(2940, -420),
(2044, 220),
$

$
(735, 1155)\}
$
(Figure~\ref{picture_1260.png}),

where $n = 9$, $\operatorname{diam(\mathcal{P})} = 1260$. Using the ``blowing up''
construction, we obtain
\begin{equation}\label{result3}
d(m, 4m + 1) \leq 1260
\end{equation}
Figure~\ref{picture_1260_R3.pdf} shows an example for $m = 3$.

However, we have to admit that both bounds which employ pyramid sets
are indeed weaker than the one
\begin{equation}
	d(m, m^2+m)\leq 17
\end{equation}
proved in~\cite{kurz2008bounds},
so the constructions of bounds based on pyramid sets tend to be only for the conceptual purpose.

\end{itemize}

\section{Final remarks}
All the discussed PIPSs were obtained through a combination of computer search and intuition of the authors;
so, the further search may lead to better bounds employing the same constructions.

There is still no general construction for a rails or scissors PIPS of arbitrary cardinality.
For rails PIPSs, we can conjecture that there exists a set of arbitrary cardinality, with 2 points on one line
and all the rest on the other;
on the other hand, it is still unknown whether there are any rails PIPSs with 4 and 5 points on the lines.

As for today, we have found pyramid PIPS of cardinality at most 9.

The source code can be obtained at https://gitlab.com/Nickkolok/ips-algo

\section{Acknowledgements}
The authors thank A.S. Usachev for his help in polishing the text of the paper.


\begin{thebibliography}{99}
{}
\bibitem{ascher2020erdos}
\emph{Ascher} \emph{K.}, \emph{Braune} \emph{L.}, \emph{Turchet} \emph{A.}
The Erdős–Ulam problem, Lang’s conjecture and uniformity // Bulletin of the London
Mathematical Society. — 2020. — Vol. 52, no. 6. — Pp. 1053–1063. — arXiv: \href
{http://arxiv.org/abs/1809.02037} {\nolinkurl {1809.02037}}.
{}\bibitem{solymosi2010question}
\emph{Solymosi} \emph{J.}, \emph{De Zeeuw} \emph{F.} On a question of Erdős
and Ulam // Discrete \& Computational Geometry. — 2010. — Vol. 43, no. 2. —
Pp. 393–401. — arXiv: \href {http://arxiv.org/abs/0806.3095} {\nolinkurl
{0806.3095}}.
{}
\bibitem{kurz2007enumeration}
\emph{Kurz} \emph{S.} Enumeration of Integral Tetrahedra // Journal of Integer
Sequences. — 2007. — Vol. 10, no. 2. — P. 3.
{}
\bibitem{anning1945integral}
\emph{Anning} \emph{N. H.}, \emph{Erdős} \emph{P.} Integral distances //
Bulletin of the American Mathematical Society. — 1945. — Vol. 51, no. 8. —
Pp. 598–600. — DOI: \href {https://doi.org/10.1090/S0002-9904-1945-08407-9}
{\nolinkurl {10.1090/S0002-9904-1945-08407-9}}.
{}
\bibitem{erdos1945integral}
\emph{Erdős} \emph{P.} Integral distances // Bulletin of the American Mathematical
Society. — 1945. — Vol. 51, no. 12. — P. 996. — DOI: \href
{https://doi.org/10.1090/S0002-9904-1945-08490-0} {\nolinkurl
{10.1090/S0002-9904-1945-08490-0}}.
{}
\bibitem{kemnitz1988punktmengen}
\emph{Kemnitz} \emph{A.} Punktmengen mit ganzzahligen Abständen. — 1988.
{}
\bibitem{kurz2005characteristic}
\emph{Kurz} \emph{S.} On the characteristic of integral point sets in {$\mathbb
{E}^m$} // Australasian Journal of Combinatorics. — 2006. — Vol. 36. —
Pp. 241–248. — arXiv: \href {http://arxiv.org/abs/math/0511704} {\nolinkurl
{math/0511704}}.
{}
\bibitem{kurz2008bounds}
\emph{Kurz} \emph{S.}, \emph{Laue} \emph{R.} Bounds for the minimum diameter
of integral point sets // Australasian Journal of Combinatorics. — 2007. — Vol. 39. —
Pp. 233–240. — arXiv: \href {http://arxiv.org/abs/0804.1296} {\nolinkurl
{0804.1296}}.
{}
\bibitem{solymosi2003note}
\emph{Solymosi} \emph{J.} Note on integral distances // Discrete \& Computational
Geometry. — 2003. — Vol. 30, no. 2. — Pp. 337–342. — DOI: \href
{https://doi.org/10.1007/s00454-003-0014-7} {\nolinkurl {10.1007/s00454-003-0014-7}}.{}
\bibitem{my-pps-linear-bound-2019}
\emph{Avdeev} \emph{N.} On existence of integral point sets and their diameter
bounds // Australasian Journal of Combinatorics. — 2020. — Vol. 77, no. 1. —
Pp. 100–116. — arXiv: \href {http://arxiv.org/abs/1906.11926} {\nolinkurl
{1906.11926}}. — URL: \url
{https://ui.adsabs.harvard.edu/abs/2019arXiv190611926A/abstract}.
{}
\bibitem{nozaki2013lower}
\emph{Nozaki} \emph{H.} Lower bounds for the minimum diameter of integral point
sets // Australasian Journal of Combinatorics. — 2013. — Vol. 56. — Pp. 139–143.
{}
\bibitem{our-vmmsh-2018}
\emph{Авдеев} \emph{Н. Н.}, \emph{Семёнов} \emph{Е. М.} Множества точек с
целочисленными расстояниями на плоскости и в евклидовом пространстве //
Математический форум (Итоги науки. Юг России). — 2018. — Pp. 217–236.
{}
\bibitem{antonov2008maximal}
\emph{Antonov} \emph{A. R.}, \emph{Kurz} \emph{S.} Maximal integral point sets
over {$\mathbb {Z}^2$} // International Journal of Computer Mathematics. —
2008. — Vol. 87, no. 12. — Pp. 2653–2676. — DOI: \href
{https://doi.org/10.1080/00207160902993636} {\nolinkurl
{10.1080/00207160902993636}}. — arXiv: \href {http://arxiv.org/abs/0804.1280}
{\nolinkurl {0804.1280}}.
{}
\bibitem{kurz2008minimum}
\emph{Kurz} \emph{S.}, \emph{Wassermann} \emph{A.} On the minimum diameter
of plane integral point sets // Ars Combinatoria. — 2011. — Vol. 101. — Pp. 265–287. —
arXiv: \href {http://arxiv.org/abs/0804.1307} {\nolinkurl {0804.1307}}.
{}
\bibitem{our-ped-2018}
\emph{Авдеев} \emph{Н. Н.} On integral point sets in special position // Некоторые
вопросы анализа, алгебры, геометрии и математического образования: материалы
международной молодежной научной школы «Актуальные направления
математического анализа и смежные вопросы». — 2018. — Vol. 8. — Pp. 5–6.
{}
\bibitem{our-pmm-2018}
\emph{Авдеев} \emph{Н. Н.} Об отыскании целоудалённых множеств
специального вида // Актуальные проблемы прикладной математики,
информатики и механики - сборник трудов Международной научной
конференции. — Научно-исследовательские публикации, 2018. — Pp. 492–498.
{}\bibitem{harborth1993upper}
\emph{Harborth} \emph{H.}, \emph{Kemnitz} \emph{A.},
\emph{Möller} \emph{M.} An upper bound for the minimum diameter of integral
point sets // Discrete \& Computational Geometry. — 1993. — Vol. 9, no. 4. —
Pp. 427–432. — DOI: \href {https://doi.org/10.1007/bf02189331} {\nolinkurl
{10.1007/bf02189331}}.
{}
\bibitem{piepmeyer1996maximum}
\emph{Piepmeyer} \emph{L.} The maximum number of odd integral distances
between points in the plane // Discrete \& Computational Geometry. — 1996. —
Vol. 16, no. 1. — Pp. 113–115. — DOI: \href {https://doi.org/10.1007/bf02711135}
{\nolinkurl {10.1007/bf02711135}}.
{}
\bibitem{bat2018number}
\emph{Bat-Ochir} \emph{G.} On the number of points with pairwise integral
distances on a circle // Discrete Applied Mathematics. — 2018. — Vol. 254. —
Pp. 17–32. — DOI: \href {https://doi.org/10.1016/j.dam.2018.07.004} {\nolinkurl
{10.1016/j.dam.2018.07.004}}.
{}
\bibitem{kreisel2008heptagon}
\emph{Kreisel} \emph{T.}, \emph{Kurz} \emph{S.} There are integral heptagons, no
three points on a line, no four on a circle // Discrete \& Computational Geometry. —
2008. — Vol. 39, no. 4. — Pp. 786–790. — DOI: \href
{https://doi.org/10.1007/s00454-007-9038-6} {\nolinkurl {10.1007/s00454-007-9038-6}}.
{}
\bibitem{kurz2013constructing}
Constructing $7$-clusters / S. Kurz [et al.] // Serdica Journal of Computing. — 2014. —
Vol. 8, no. 1. — Pp. 47–70. — arXiv: \href {http://arxiv.org/abs/1312.2318} {\nolinkurl
{1312.2318}}.
{}
\bibitem{avdeev2021particular}
\emph{Avdeev} \emph{N.}, \emph{Momot} \emph{E.},
\emph{Zvolinsky} \emph{A.} On particular examples of planar integral point sets and
their classification. —. — arXiv: \href {http://arxiv.org/abs/2102.12462 [math.CO]}
{\nolinkurl {2102.12462 [math.CO]}}.
{}
\bibitem{my-semi-general-5-4-bound-2019}
\emph{Avdeev} \emph{N.} On diameter bounds for planar integral point sets in
semi-general position. — 2019. — arXiv: \href {http://arxiv.org/abs/1907.09331}
{\nolinkurl {1907.09331}}. — URL: \url
{https://ui.adsabs.harvard.edu/abs/2019arXiv190709331A/abstract}.
{}\bibitem{guy2013unsolved}
\emph{Guy} \emph{R.} Unsolved problems in number theory. Vol. 1. — Springer
Science \& Business Media, 2013. — DOI: \href
{https://doi.org/10.1007/978-1-4757-1738-9} {\nolinkurl {10.1007/978-1-4757-1738-9}}.
\end{thebibliography}

\end{document}